\documentclass[reqno]{amsart}

\usepackage{amsmath}
\usepackage{amsfonts}
\usepackage{amssymb}
\usepackage{array}
\usepackage{tabularx}
\usepackage{bbm}
\usepackage{url}
\usepackage{graphicx}
\usepackage{psfrag}
\usepackage{booktabs}
\usepackage[labelfont={footnotesize,bf},textfont=footnotesize]{caption}
\usepackage[labelfont=footnotesize,textfont=footnotesize]{subfig}
\usepackage[T1]{fontenc}
\usepackage[noend]{myalgorithmic}
\usepackage{algorithm}
\usepackage{mathtools}

\usepackage{pdfsync}

\graphicspath{{Pictures/}}

\newtheorem{theorem}{Theorem}
\newtheorem{remark}{Remark}
\newtheorem{lemma}{Lemma}
\newtheorem{corollary}{Corollary}
\newtheorem{proposition}{Proposition}

\DeclareMathOperator{\convop}{conv}
\DeclareMathOperator{\bigoop}{O}
\DeclareMathOperator{\affop}{aff}

\DeclareMathOperator{\cubeop}{C}

\DeclareMathOperator{\fixop}{Fix}
\DeclareMathOperator{\red}{red}
\DeclareMathOperator{\green}{green}
\DeclareMathOperator{\white}{white}
\DeclareMathOperator{\orbiop}{O}
\DeclareMathOperator{\rowop}{row}

\DeclareMathOperator{\polytime}{P}
\DeclareMathOperator{\NP}{NP}
\DeclareMathOperator{\coNP}{coNP}

\newcommand{\ints}[1]{[{#1}]}
\newcommand{\setdef}[2]{\{{#1}\,:\,{#2}\}}

\newcommand{\R}{\mathbbm{R}}
\newcommand{\Q}{\mathbbm{Q}}
\newcommand{\N}{\mathbbm{N}}

\newcommand{\bigo}[1]{\bigoop({#1})}
\newcommand{\aff}[1]{\affop({#1})}
\newcommand{\card}[1]{\lvert {#1} \rvert}
\newcommand{\define}{\coloneqq}

\providecommand*{\dotcup}{\ensuremath{\stackrel{\cdot}{\cup}}}

\providecommand*{\unitvec}[1]{\ensuremath{\mathbbm{e}_{#1}}}

\newcommand{\row}[1]{\rowop_{#1}}

\newcommand{\cube}[1]{\ensuremath{\cubeop^{#1}}}
\newcommand{\fixface}[2]{\fixop_{#2}({#1})}
\newcommand{\fixfaceaff}[3]{\fixop^{#3}_{#2}({#1})}
\newcommand{\orbipart}[2]{\ensuremath{\orbiop^{=}_{#1,#2}}}
\newcommand{\orbipack}[2]{\ensuremath{\orbiop^{\le}_{#1,#2}}}
\newcommand{\orbicov}[2]{\ensuremath{\orbiop^{\ge}_{#1,#2}}}

\newcommand{\orbipartinds}[2]{\ensuremath{\mathcal{I}_{{#1},{#2}}}}
\newcommand{\diagcol}[2]{\langle{#1},{#2}\rangle}
\newcommand{\syssci}{\mathcal{S}_{\text{SCI}}}
\newcommand{\sysci}{\mathcal{S}_{\text{CI}}}

\newenvironment{myitemize}{%
\begin{list}{$\circ$}%
{\setlength{\topsep}{0.5ex}%
\setlength{\partopsep}{0mm}%
\setlength{\parskip}{0ex}%
\setlength{\parsep}{0mm}%
\setlength{\itemsep}{0ex}%
\setlength{\labelwidth}{4mm}%
\setlength{\leftmargin}{0ex}%
\addtolength{\leftmargin}{\labelwidth}%
\addtolength{\leftmargin}{\labelsep}%
\setlength{\itemindent}{0mm}}}%
{\end{list}}

\makeatletter
\def\@secnumfont{\bfseries}
\renewcommand\section{\@startsection {section}{1}{\z@}%
                                   {-3.5ex \@plus -1ex \@minus -.2ex}%
                                   {2.3ex \@plus.2ex}%
                                   {\normalfont\large\bfseries}}
\renewcommand\subsection{\@startsection{subsection}{2}{\z@}%
                                     {-3.25ex\@plus -1ex \@minus -.2ex}%
                                     {1.5ex \@plus .2ex}%
                                     {\normalfont\normalsize\bfseries}}
\renewcommand\subsubsection{\@startsection{subsubsection}{3}{\z@}%
                                     {-3.25ex\@plus -1ex \@minus -.2ex}%
                                     {1.5ex \@plus .2ex}%
                                     {\normalfont\normalsize\bfseries}}
\makeatother

\begin{document}

\title{Orbitopal Fixing}
\thanks{This work has been partially supported by the DFG Research Center
  \textsc{Matheon} in Berlin.}

\author{Volker Kaibel}
\thanks{During part of the research presented in this paper the first author
  was a visiting professor at TU Berlin}
\address{
  Otto-von-Guericke Universit\"at Magdeburg \\
  Fakult\"at f\"ur Mathematik \\
  Universit\"atsplatz 2 \\
  39106 Magdeburg \\
  Germany
}
\email{kaibel@ovgu.de}

\author{Matthias Peinhardt}
\address{
  Otto-von-Guericke Universit\"at Magdeburg \\
  Fakult\"at f\"ur Mathematik \\
  Universit\"atsplatz 2 \\
  39106 Magdeburg \\
  Germany
}
\email{peinhard@ovgu.de}

\author{Marc\,E. Pfetsch}
\address{%
  TU Braunschweig \\
  Institute for Mathematical Optimization \\
  Pockelsstr.\ 14 \\
  38106 Braunschweig \\
  Germany
}
\email{m.pfetsch@tu-bs.de}

\begin{abstract}
  The topic of this paper are integer programming models in which a subset
  of 0/1-variables encode a partitioning of a set of objects into disjoint
  subsets.  Such models can be surprisingly hard to solve by branch-and-cut
  algorithms if the order of the subsets of the partition is irrelevant, since
  this kind of symmetry unnecessarily blows up the search tree.

  We present a general tool, called orbitopal fixing, for enhancing the
  capabilities of branch-and-cut algorithms in solving such symmetric
  integer programming models. We devise a linear time algorithm that,
  applied at each node of the search tree, removes redundant parts
  of the tree produced by the above mentioned symmetry. The method relies
  on certain polyhedra, called orbitopes, which have been introduced
  in~\cite{KaPf08}. It does, however, not explicitly add inequalities to the model.
  Instead, it uses certain fixing rules for variables.
  We demonstrate the computational power of
  orbitopal fixing at the example of a graph partitioning problem.
\end{abstract}

\maketitle

\section{Introduction}
\label{sec:introduction}

Being welcome in most other contexts, symmetry causes severe trouble in the
solution of many integer programming (IP) models. This paper describes a
method to enhance the capabilities of branch-and-cut algorithms with
respect to hand\-ling symmetric models of a certain kind that frequently
occurs in practice.

We illustrate this kind of symmetry by the example of a graph partitioning
problem (another notorious example is the vertex coloring problem).  Here,
one is given an undirected graph $G=(V,E)$ with non-negative edge weights
$w\in\Q_{\ge 0}^E$ and an integer $q\ge 2$. The task is to partition~$V$
into~$q$ disjoint subsets such that the sum of all weights of edges
connecting nodes in the same subset is minimized; thus, this problem is
equivalent to maximizing the weights of the edges in a $q$-cut.

A straight-forward IP model for this graph partitioning problem arises by
introducing 0/1-variables $x_{ij}$ for all $i\in\ints{p} \define
\{1,\dots,p\}$ and $j\in\ints{q}$ that indicate whether node~$i$ is
contained in subset~$j$ (where we assume $V=\ints{p}$).  In order to model
the objective function, we furthermore need 0/1-variables $y_{ik}$, for all
edges $\{i,k\}\in E$, indicating whether nodes~$i$ and~$k$ are contained in
the same subset. This yields the following IP-model (see, e.g.,
\cite{Eis01}):

\begin{equation}\label{eq:kpart}
\begin{alignedat}{3}
&\min        & \sum_{\{i,k\}\in E}w_{ik}\,y_{ik} &             && \\
&\text{s.t.} & \sum_{j=1}^q x_{ij}    & = 1    && \qquad\text{ for all }i\in \ints{p} \\
&            &  x_{ij}+x_{kj}-y_{ik}  & \le 1  && \qquad\text{ for all }\{i,k\} \in E,\, j\in\ints{q}\\
&            & x_{ij}                 & \in \{0,1\} && \qquad\text{ for all }i\in \ints{p},\,j\in\ints{q}\\
&            & y_{ik}                 & \in \{0,1\} && \qquad\text{ for all }\{i,k\}\in E.
\end{alignedat}
\end{equation}

The $x$-variables describe a 0/1-matrix of size $p\times q$ with exactly
one $1$-entry per row. They encode the assignment of the nodes to the
subsets of the partition. The methods that we discuss in this paper do only
rely on this structure and thus can be applied to many other models as
well. We use the example of the graph partitioning problem as a prototype
application and report on computational experiments for this application in
Sect.~\ref{sec:compu}.

Graph partitioning problems are discussed, for example, in
\cite{ChRa93,ChRa95,Eis01}. They arise, for instance, as  relaxations of
frequency assignment problems in mobile telecommunication networks,
see~\cite{Eis01}. The maximization version (min $k$-cut)
of the graph
partitioning problem is relevant as well~\cite{FaReWo94,KoGlAlWa05}. Also
capacity bounds on the subsets of the partition (which can easily be
incorporated into the model) are of interest, in particular for the graph
equipartition problem \cite{FeMaSoWeWo96,FeMaSoWeWo98,MeTr98,So05}. For
the closely related clique partitioning problem see~\cite{GrWa89,GrWa90}.
Semidefinite relaxations and solution approaches are discussed
in~\cite{Eis02,GahAL09}.

As it is given above, the model is unnecessarily difficult for
state-of-the-art IP solvers. Even solving small instances requires enormous
efforts (see Sect.~\ref{sec:compu}).  One reason is that every feasible
solution $(x,y)$ to this model can be turned into~$q!$ different ones by
permuting the columns of~$x$ (viewed as a 0/1-matrix) in an arbitrary way,
thereby not changing the structure of the solution (in particular: its
objective function value).  Phrased differently, the symmetric group of all
permutations of the set~$\ints{q}$ operates on the solutions by permuting
the columns of the $x$-variables in such a way that the objective function
remains constant along each orbit.  Therefore, when solving the model by a
branch-and-cut algorithm, basically the same work will be done in the tree
at many places.  Thus, there should be potential for reducing the running
times significantly by exploiting this symmetry.  A more subtle second point
is that interior points of the convex hulls of the individual orbits are
responsible for quite weak linear programming (LP) bounds. We will,
however, not address this second point in this paper.

In order to remove symmetry, the above model for the graph partitioning
problem is often replaced by models containing only edge variables, see,
e.g.~\cite{FeMaSoWeWo96}. For this, however, the underlying graph
has to be complete, which might introduce many unnecessary variables.
Moreover, formulation~(\ref{eq:kpart}) is sometimes favorable, e.g., if
node-weighted capacity constraints should be incorporated.

One way to deal with symmetry is to restrict the feasible region in each of
the orbits to a single representative, e.g., to the lexicographically
maximal (with respect to the row-by-row ordering of the $x$-components)
element in the orbit. In fact, this can be done by adding inequalities to
the model that enforce the columns of $x$ to be sorted in a
lexicographically decreasing way. This can be achieved by $\bigo{pq}$ many
\emph{column inequalities}. In~\cite{KaPf08} even a complete (and
irredundant) linear description of the convex hull of all 0/1-matrices of
size $p\times q$ with exactly one $1$-entry per row and lexicographically
decreasing columns is derived; a shorter proof of this completeness result
appears in~\cite{FaeK09}. The corresponding polytopes are called
\emph{partitioning orbitopes}. A similar result can be proved for the case
of \emph{packing orbitopes}, in which there is at most one $1$-entry per
row. The descriptions basically consist of an exponentially large
super class of column inequalities, called \emph{shifted column
  inequalities}, for which there is a linear time separation algorithm
available. We recall some of these results in Sect.~\ref{sec:orbitopes}.

Incorporating the inequalities from the orbitope description into the IP
model removes symmetry. At each node of the branch-and-cut tree this
ensures that the corresponding IP is infeasible as soon as there is no
representative in the subtree rooted at that node. In fact, already the
column inequalities are sufficient for this purpose.

In this paper, we investigate a way to utilize these inequalities (or the
orbitope that they describe) without explicitly adding any of the
inequalities to the models. The reason for doing this is the unpleasant
effect that adding (shifted) column inequalities to the models might result
in more difficult LP relaxations. One way of avoiding the addition of these
inequalities to the LPs is to derive logical implications instead: If we
are working in a branch-and-cut node at which the $x$-variables
corresponding to index subsets~$I_0$ and~$I_1$ are fixed to zero and one,
respectively, then there might be a (shifted) column inequality yielding
implications for all representatives in the subtree rooted at the current
node. For instance, it might be (and this is easy to check for a given
inequality) that for some $(i^{\star},j^{\star})\not\in I_0\cup I_1$ we
have $x_{i^{\star} j^{\star}}=0$ for all 0/1-points~$x$ with $x_{ij}=0$
($(i,j)\in I_0$) and $x_{ij}=1$ ($(i,j)\in I_1$) that satisfy the
inequality. In this case, $x_{i^{\star} j^{\star}}$ can be fixed to zero
for the whole subtree rooted at the current node, enlarging~$I_0$.
Similarly, also fixings of variables to 1 might be possible.  We call the
iterated process of searching for such additional fixings \emph{sequential
  fixing} with (shifted) column inequalities.

Let us mention at this point that deviating from parts of the literature,
we do not distinguish between ``fixing'' and ``setting'' of variables in
this paper.

Sequential fixing with (shifted) column inequalities is a special case of
constraint propagation, which is well known from constraint logic
programming, see~\cite{Apt03,Hen89,MarS98} for an overview. Modern IP
solvers like SCIP~\cite{Ach07} use such strategies also in the node
preprocessing during the branch-and-cut algorithm. With orbitopes, however,
we can aim at something better: Consider a branch-and-cut node identified
by fixing the variables corresponding to sets~$I_0$ and~$I_1$ to zero and
one, respectively. Denote by $W(I_0,I_1)$ the set of all vertices~$x$ of
the orbitope with $x_{ij}=0$ for all $(i,j)\in I_0$ and $x_{ij}=1$ for all
$(i,j)\in I_1$. We define the sets $I^{\star}_0$ and $I^{\star}_1$ of all
indices $(i^{\star},j^{\star})$ of variables, for which all~$x$ in
$W(I_0,I_1)$ satisfy $x_{i^{\star}j^{\star}}=0$ and
$x_{i^{\star}j^{\star}}=1$, respectively. We call the respective fixing of
the variables corresponding to $I^{\star}_0$ and~$I^{\star}_1$
\emph{simultaneous fixing}.  Simultaneous fixing is always at least as
strong as sequential fixing.

Investigations of sequential and simultaneous fixing for orbitopes are the
central topic of the paper.  The main contributions and results are the
following:
\begin{myitemize}
\item We present a linear time algorithm for \emph{orbitopal fixing}, i.e.,
  for solving the problem to compute simultaneous fixings for partitioning
  orbitopes (Theorem~\ref{thm:algorbifix}) and packing orbitopes
  (Corollary~\ref{cor:PackingOrbitopes}).
\item In contrast to this, we prove that orbitopal fixing for
  \emph{covering orbitopes} (the convex hulls of all lexicographically
  maximal 0/1-matrices with \emph{at least} one $1$-entry in every row) is
  $\NP$-hard (Theorem~\ref{thm:orbicov}).
\item We show that, for general 0/1-polytopes, sequential fixing, even with
  complete and irredundant linear descriptions, is weaker than simultaneous
  fixing (Theorem~\ref{thm:fix}). For the case of partitioning orbitopes,
  we clarify the relationships between different versions of sequential
  fixing with (shifted) column inequalities, where (despite the situation
  for general 0/1-polytopes) the strongest one is as strong as orbitopal
  fixing (Theorem~\ref{thm:seqfixOrbi}).
\item We report on computer experiments (Sect.~\ref{sec:compu}) with the
  graph partitioning problem described above, showing that orbitopal fixing
  leads to significant performance improvements for branch-and-cut
  algorithms.
\end{myitemize}
\smallskip

This paper extends the one that appeared in the proceedings of IPCO~XII~\cite{KaPePf07}. It
contains the following additional material: a proof for the
second part of Theorem~\ref{thm:seqfixOrbi}, the above mentioned results
for packing and covering orbitopes (Sect.~\ref{sec:packcov}), a comparison
to the related approaches of Margot~\cite{Mar02,Mar03b,Mar07} and Linderoth et al.~\cite{LiOsRoSm06,OstLRS08} for the orbitope
case (Sect.~\ref{sec:Margot}), and, finally,
computational results for a significantly improved version of our graph
partitioning code (Sect.~\ref{sec:compu}).

While our methods are based on lexicographically maximal choices of
representatives from the orbits, a more general approach admitting
orderings defined by arbitrary linear functions was introduced by Friedman,
see~\cite{Fr07}. There are also a number of approaches for symmetry
handling available from the constraint logic programming literature, see,
e.g., \cite{FahSS01b,Pug05,SelH05}. Their general idea is similar to the
above mentioned approaches by Margot and Linderoth et al. During the
traversal of the tree, different techniques are used to avoid the
processing of (some) symmetric parts of the tree.
For an excellent survey of methods for symmetry breaking in integer
programming we refer to~\cite{Mar10}.

\section{Orbitopes}\label{sec:orbitopes}

Throughout the paper, let~$p$ and~$q$ be integers with $p\ge q\ge 2$. The
\emph{partitioning/packing/covering orbitope}
$\orbipart{p}{q}$/$\orbipack{p}{q}$/$\orbicov{p}{q}$ is the convex hull of
all 0/1-matrices $x\in\{0,1\}^{\ints{p}\times\ints{q}}$ with exactly/at
most/at least one $1$-entry per row, whose columns are in non-increasing
lexicographical order, i.e., they satisfy
\begin{equation}\label{eq:lexorder}
\sum_{i=1}^p 2^{p-i} x_{ij} \ge \sum_{i=1}^p 2^{p-i}x_{i,j+1}
\end{equation}
for all $j\in\ints{q-1}$.

We will mainly be concerned with partitioning orbitopes~$\orbipart{p}{q}$.
An exception is Sect.~\ref{sec:packcov}, in which we will show that the
linear time method for orbitopal fixing of Sect.~\ref{sec:orbifix} below
can  easily be carried over to packing orbitopes~$\smash{\orbipack{p}{q}}$, while
there is no polynomial time method for orbitopal fixing for covering
orbitopes~$\smash{\orbicov{p}{q}}$, unless $\polytime=\NP$.

Let the symmetric group of size~$q$ act on
$\{0,1\}^{\ints{p}\times\ints{q}}$ via permutation of the columns. Then the
vertices of $\orbipart{p}{q}$ are exactly the lexicographically maximal
matrices with exactly one $1$-entry per row in the orbits under the
symmetric group action; the lexicographic order is defined as
in~\eqref{eq:lexorder}.

 As these vertices have $x_{ij}=0$ for all $(i,j)$ with $i<j$, we drop
these components and consider $\orbipart{p}{q}$ as a subset of the space
$\R^{\orbipartinds{p}{q}}$ with $\orbipartinds{p}{q} \define \setdef{(i,j)
  \in \{0,1\}^{\ints{p}\times\ints{q}}}{i \geq j}$. Thus, we consider
matrices, in which the $i$-th row has $q(i) \define \min\{i,q\}$
components.

\begin{figure}[tb]
  \centering
  \newcommand{\mystyle}[1]{\footnotesize {#1}}
  \psfrag{i}{\mystyle{$i$}}
  \psfrag{j}{\mystyle{$j$}}
  \psfrag{eta}{\mystyle{$\eta$}}
  \subfloat[][]{\includegraphics[width=.15\textwidth]{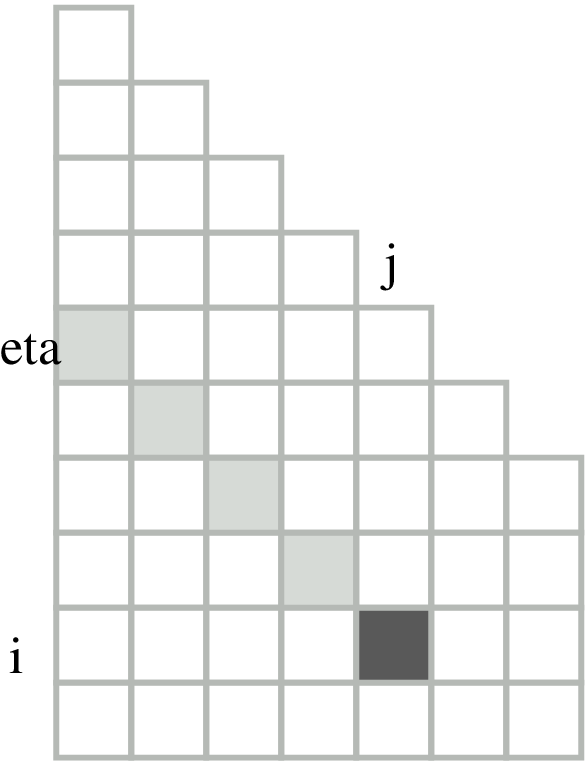}\label{sfl:coord}}\hspace{8ex}
  \subfloat[][]{\includegraphics[width=.15\textwidth]{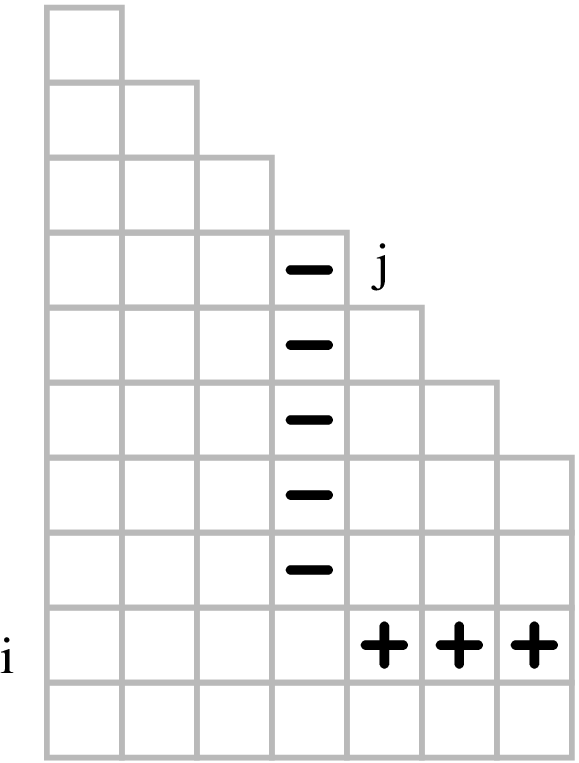}\label{sfl:scipart1}}\hspace{8ex}
  \subfloat[][]{\includegraphics[width=.15\textwidth]{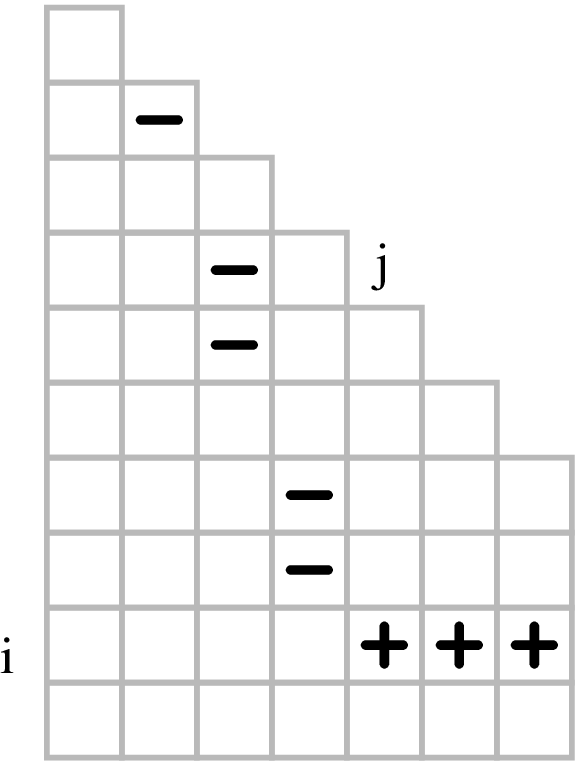}\label{sfl:scipart2}}\hspace{8ex}
  \subfloat[][]{\includegraphics[width=.15\textwidth]{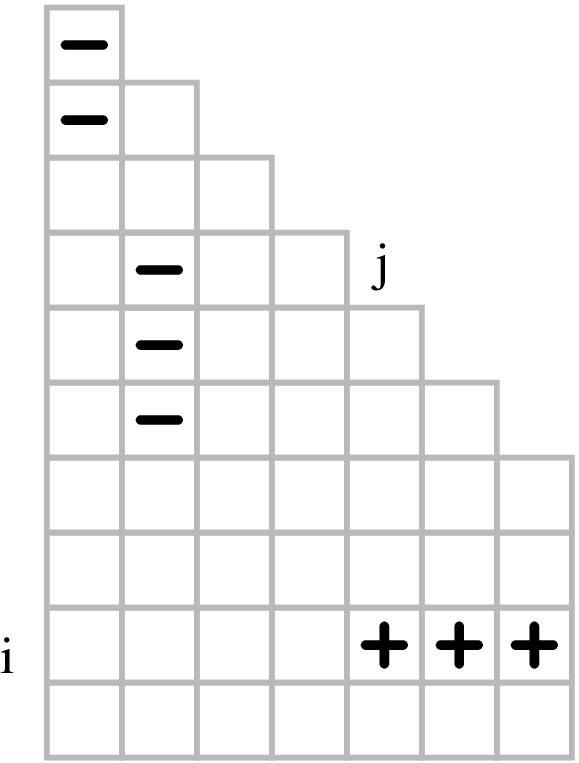}\label{sfl:scipart3}}
  \caption[]{\subref{sfl:coord} Example for coordinates $(9,5) =
    \diagcol{5}{5}$. \subref{sfl:scipart1}, \subref{sfl:scipart2},
    \subref{sfl:scipart3} Three shifted column inequalities,
    \subref{sfl:scipart1} being a column inequality}
  \label{fig:scispart}
\end{figure}

The main result in~\cite{KaPf08} is a complete linear description of
$\orbipart{p}{q}$. In order to describe the result, it will be convenient
to address the elements in $\orbipartinds{p}{q}$ via a different ``system
of coordinates'': For $j \in \ints{q}$ and $1 \leq \eta \leq p-j+1$, define
$\diagcol{\eta}{j} \define (j+\eta-1,j)$. Thus (as before) $i$ and~$j$
denote the row and the column, respectively, while~$\eta$ is the index of
the diagonal (counted from above) containing the respective element; see
Figure~\ref{fig:scispart}~\subref{sfl:coord} for an example.

A set $S = \{\diagcol{1}{c_1}, \diagcol{2}{c_2}, \dots,
\diagcol{\eta}{c_{\eta}}\} \subset \orbipartinds{p}{q}$ with $c_1 \leq c_2
\leq \dots \leq c_{\eta}$ and $\eta \geq 1$ is called a \emph{shifted
  column}.  For $(i,j) = \diagcol{\eta}{j} \in \orbipartinds{p}{q}$, a
shifted column~$S$ as above with $c_{\eta}<j$, and the set $B = \{(i,j),
(i,j+1), \dots, (i,q(i))\}$, we call $x(B) - x(S) \leq 0$ a \emph{shifted
  column inequality}. The set~$B$ is called its \emph{bar}. In case of
$c_1=\dots=c_{\eta}=j-1$ the shifted column inequality is called a
\emph{column inequality}.  See Figure~\ref{fig:scispart} for examples.

Finally, a bit more notation is needed. For each $i \in \ints{p}$, we
define $ \row{i} \define \setdef{(i,j)}{j\in\ints{q(i)}}$. For $A \subset
\orbipartinds{p}{q}$ and $x \in \R^{\orbipartinds{p}{q}}$, we denote
by~$x(A)$ the sum $\sum_{(i,j)\in A}x_{ij}$.

\begin{theorem}[see \cite{KaPf08}]\label{thm:OrbitpartDesc}
  The orbitope $\orbipart{p}{q}$ is completely described by the
  non-negativity constraints $x_{ij}\ge 0$, the row-sum equations
  $x(\row{i})=1$, and the shifted column inequalities.
\end{theorem}

In fact, in~\cite{KaPf08} it is also shown that, up to a few exceptions,
the inequalities in this description define facets of $\orbipart{p}{q}$.
Furthermore, a linear time separation algorithm for the exponentially large
class of shifted column inequalities is given. For a compact extended
formulation of $\orbipart{p}{q}$ that also leads to a simplified proof of
Theorem~\ref{thm:OrbitpartDesc}, see~\cite{FaeK09}.

\section{The Geometry of Fixing Variables}
\label{sec:geomfix}

In this section, we deal with general 0/1-integer programs and, in
particular, their associated polytopes. We will define some basic
terminology used later in the special treatment of orbitopes, and we are
going to shed some light on the geometric situation of fixing variables.

For some positive integer $d$, we denote by
\[
\cube{d} = \setdef{x\in\R^d}{0\le x_i \le 1 \text{ for all }i \in \ints{d}}
\]
the 0/1-cube, where $\ints{d}$ is the corresponding set of indices of
variables. For two disjoint subsets $I_0, I_1 \subseteq\ints{d}$ (hence,
$I_0 \cap I_1 = \varnothing$) we call
\[
\setdef{x\in\cube{d}}{x_i=0 \text{ for all }i\in I_0,\
  x_i=1 \text{ for all }i\in I_1}
\]
the \emph{face of~$\cube{d}$ defined by $(I_0,I_1)$}. All nonempty faces of
$\cube{d}$ are of this type.

For a polytope $P\subseteq\cube{d}$ and for a face~$F$ of $\cube{d}$
defined by $(I_0,I_1)$, we denote by $\fixface{P}{F}$ the smallest face
of~$\cube{d}$ that contains $P\cap F\cap\{0,1\}^d$ (i.e., $\fixface{P}{F}$
is the intersection of all faces of~$\cube{d}$ that contain $P\cap
F\cap\{0,1\}^d$). If $\fixface{P}{F}$ is the nonempty cube face defined by
$(I^{\star}_0,I^{\star}_1)$, then $I^{\star}_0$ and $I^{\star}_1$ consist
of all $i\in\ints{d}$ for which $x_i=0$ and $x_i=1$, respectively, holds
for all $x\in P\cap F\cap\{0,1\}^d$. In particular, we have $I_0\subseteq
I^{\star}_0$ and $I_1\subseteq I^{\star}_1$, or $\fixface{P}{F} =
\varnothing$.  Thus, if $I_0$ and $I_1$ are the indices of the variables
fixed to zero and one, respectively, in the current branch-and-cut node
(with respect to an IP with feasible points $P\cap\{0,1\}^d$), the node can
either be pruned, or the sets $I^{\star}_0$ and $I^{\star}_1$ yield the
maximal sets of variables that can be fixed to zero and one, respectively,
for the whole subtree rooted at this node. Unless
$\fixface{P}{F}=\varnothing$, we call $(I^{\star}_0,I^{\star}_1)$ the
\emph{fixing of~$P$ at $(I_0,I_1)$}.  Similarly, we call $\fixface{P}{F}$
the \emph{fixing} of~$P$ at~$F$.

\begin{remark}\label{rem:fixmono}
  If $P,P'\subseteq\cube{d}$ are two polytopes with $P\subseteq P'$ and~$F$
  and~$F'$ are two faces of~$\cube{d}$ with $F\subseteq F'$, then
  $\fixface{P}{F}\subseteq\fixface{P'}{F'}$ holds.
\end{remark}

In general, it is not clear how to compute fixings efficiently. Indeed,
computing the fixing of~$P$ at $(\varnothing,\varnothing)$ includes
deciding whether~$P\cap\{0,1\}^d=\varnothing$, which, of course, is $\NP$-hard
in general. On the other hand, the following holds.

\begin{lemma}\label{lemma:PolyTimeFixing}
  If one can optimize a linear function over~$P \cap \{0,1\}^d$ in
  polynomial time, the fixing $(I_0^\star, I_1^\star)$ at $(I_0, I_1)$ can
  be computed in polynomial time.
\end{lemma}

\begin{proof}
  Let $c \in \R^d$ be the objective function vector defined by
  \[
  c_i = \begin{cases}
    \phantom{-}1 & \text{if } i \in I_1\\
    -1 & \text{if } i \in I_0\\
    \phantom{-}0 & \text{otherwise}
  \end{cases}
  \qquad
  \text{for all } i \in [d].
  \]
  For each $i^{\star} \in [d] \setminus (I_0 \cup I_1)$ we have
  \begin{equation*}
    \max\setdef{(c+\unitvec{i^{\star}})^{T}x}{x\in P\cap\{0,1\}^d}\ <\ \card{I_1}+1
  \end{equation*}
  (where $\unitvec{i}$ is the $i$th unit vector) if and only if $i^{\star}\in I^{\star}_0$, and
  \begin{equation*}
    \max\setdef{(c-\unitvec{i^{\star}})^{T}x}{x\in P\cap\{0,1\}^d}\ <\ \card{I_1}
  \end{equation*}
  if and only if $i^{\star}\in I^{\star}_1$. Thus, we can compute
  $I^{\star}_0$ and $I^{\star}_1$ by solving $2(d-\card{I_0}-\card{I_1})$
  many linear optimization problems over $P\cap\{0,1\}^d$.
\end{proof}

Note that the reverse to the implication stated in
Lemma~\ref{lemma:PolyTimeFixing} does not hold, in general. This can, e.g.,
be seen at the example of $0/1$-knapsack problems with $P = \setdef{x \in
  \R^d}{\sum_{i=1}^d a_i x_i \leq b}$ (with $a_1,\dots, a_d\ge 0$). For
every $(I_0, I_1)$ the fixing $(I_0^\star, I_1^\star)$ can be computed in
linear time:
\[
I_0^\star = I_0 \dotcup \setdef{i^{\star} \notin I_0 \cup
  I_1}{a_{i^{\star}} + \sum_{j \in I_1} a_j > b},\quad I_1^\star = I_1.
\]
In contrast, the optimization problem over $P \cap \{0,1\}^d$ is
$\NP$-hard.

If the linear optimization problem over $P\cap\{0,1\}^d$ cannot be solved
efficiently, one can still try to compute (hopefully large) subsets of
$I^{\star}_0$ and $I^{\star}_1$ by considering relaxations of~$P$. In case
of an IP that is based on an intersection with an orbitope, one might use
the orbitope as such a relaxation. We will deal with the fixing problem for
partitioning orbitopes in Sect.~\ref{sec:orbifix} (and for packing and
covering orbitopes in Sect.~\ref{sec:packcov}). Since the optimization
problem for partitioning and packing orbitopes can be solved in polynomial
time (see~\cite{KaPf08}), by Lemma~\ref{lemma:PolyTimeFixing}, the
corresponding fixing problems can be solved in polynomial time as well.
However, we will even describe linear time algorithms for these cases.

If~$P$ is given via an inequality description, one possibility is to use
the knapsack relaxations obtained from single inequalities among the
description. For each of these relaxations, the fixing can easily be
computed. If the inequality system describing~$P$ is exponentially large,
and the inequalities are only accessible via a separation routine, it might
in some cases nevertheless be possible to decide efficiently whether any of
the exponentially many knapsack relaxations allows to fix some variable
(see Sect.~\ref{subsec:fixSCI}).

Suppose, $P = \setdef{x\in\cube{d}}{Ax\leq b}$ and $P_r = \setdef{x \in
  \cube{d}}{ a_r^T x \leq b_r}$ is the knapsack relaxation of~$P$ for the
$r$th-row $a_r^T x \leq b_r$ of~$A x \leq b$, where $r = 1, \dots, m$.
Let~$F$ be some face of~$\cube{d}$. The face~$G$ of~$\cube{d}$ obtained by
setting $G \define F$ and then iteratively replacing $G$ by
$\fixface{P_r}{G}$ as long as there is some~$r\in\ints{m}$ with
$\fixface{P_r}{G}\subsetneq G$, is denoted by~$\fixface{Ax\le b}{F}$. Note
that the outcome of this procedure is independent of the choices made
for~$r$, due to Remark~\ref{rem:fixmono}. We call the pair
$(\tilde{I}_0,\tilde{I}_1)$ defining the cube face $\fixface{Ax\le b}{F}$
(unless this face is empty) the \emph{sequential fixing of~$Ax\le b$ at
  $(I_0,I_1)$}. In the context of sequential fixing we often refer to (the
computation of) $\fixface{P}{F}$ as \emph{simultaneous fixing}.

Due to Remark~\ref{rem:fixmono}, it is clear that
$\fixface{P}{F}\subseteq\fixface{Ax\le b}{F}$ holds.

\begin{theorem}\label{thm:fix}
  In general, even for a system of facet-defining inequalities describing a
  full-dimensional 0/1-polytope, sequential fixing is weaker than simultaneous fixing.
\end{theorem}

\begin{proof}
  The following example shows this. Let $P\subset\cube{4}$ be the
  four-dimensional polytope defined by the trivial inequalities $x_i\ge 0$
  for $i\in\{1,2,3\}$, $x_i \leq 1$ for $i\in\{1,2,4\}$, the inequality $
  -x_1 + x_2 + x_3 - x_4 \leq 0 $ and $x_1 - x_2 + x_3 - x_4 \leq 0.  $
  Let~$F$ be the cube face defined by $(\{4\},\varnothing)$. Then,
  sequential fixing does not fix any further variable, although
  simultaneous fixing yields $I^{\star}_0=\{3,4\}$ (and
  $I^{\star}_1=\varnothing$). Note that~$P$ has only 0/1-vertices, and all
  inequalities are facet defining ($x_4\ge 0$ and $x_3\le 1$ are
  implied).
\end{proof}

\section{Fixing Variables for Partitioning Orbitopes}
\label{sec:orbifix}

For this section, suppose that $I_0, I_1 \subseteq \orbipartinds{p}{q}$ are
subsets of indices of partitioning orbitope variables with the following properties:
\begin{myitemize}
\item[(P1)] $\card{I_0 \cap \row{i}} \leq q(i)-1$ for all $i \in \ints{p}$.
\item[(P2)] For all $(i,j) \in I_1$, we have $(i,\ell) \in I_0$ for all
  $\ell\in\ints{q(i)}\setminus \{j\}$.
\end{myitemize}
In particular, P1 and P2 imply that $I_0 \cap I_1 = \varnothing$. Let~$F$
be the face of the 0/1-cube $\cube{\orbipartinds{p}{q}}$ defined by
$(I_0,I_1)$.  Note that if P1 is not fulfilled, then $\orbipart{p}{q}\cap
F=\varnothing$.  The following statement follows immediately from
Property~P2.

\begin{remark}\label{rem:props}
  If a vertex~$x$ of $\orbipart{p}{q}$ satisfies $x_{ij}=0$ for all
  $(i,j)\in I_0$, then $x\in F$.
\end{remark}

We assume that the face $\fixface{\orbipart{p}{q}}{F}$ is defined by
$(I^{\star}_0,I^{\star}_1)$, if $\fixface{\orbipart{p}{q}}{F}$ is not empty.
\emph{Orbitopal fixing} (for partitioning orbitopes) is the problem to
compute the simultaneous fixing $(I^{\star}_0,I^{\star}_1)$ from
$(I_0,I_1)$, or determine that $\fixface{\orbipart{p}{q}}{F} = \varnothing$.

\begin{remark}\label{rem:orbiFix}
  If $\fixface{\orbipart{p}{q}}{F}\neq\varnothing$, it is enough to
  determine $I^{\star}_0$, as we have $(i,j)\in I^{\star}_1$ if and only if
  $(i,\ell)\in I^{\star}_0$ holds for for all $\ell\in\ints{q(i)}\setminus
  \{j\}$.
\end{remark}

\subsection{Intersection of Partitioning Orbitopes with Cube Faces}
\label{subsec:walls}

We start by deriving some structural results on partitioning orbitopes that
are crucial in our context.  Since
$\orbipart{p}{q}\subset\cube{\orbipartinds{p}{q}}$ is a 0/1-polytope (i.e.,
it is integral), we have $ \convop(\orbipart{p}{q}\cap
F\cap\{0,1\}^{\orbipartinds{p}{q}}) = \orbipart{p}{q}\cap F $.  Thus,
$\fixface{\orbipart{p}{q}}{F}$ is the smallest cube face that contains the
face $\orbipart{p}{q}\cap F$ of the orbitope~$\orbipart{p}{q}$.

Let us, for $i\in\ints{p}$, define values $\alpha_i \define
\alpha_i(I_0)\in\ints{q(i)}$ recursively by setting $\alpha_1\define 1$
and, for all $i\in\ints{p}$ with $i\ge 2$,
\[
\alpha_i\ \define\
\begin{cases}
\alpha_{i-1} & \text{if }\alpha_{i-1}=q(i)\text{ or }(i,\alpha_{i-1}+1)\in I_0\\
\alpha_{i-1}+1 & \text{otherwise}.
\end{cases}
\]
The set of all indices of rows, in which the $\alpha$-value increases, is
denoted by
\[
\Gamma(I_0)\ \define\ \setdef{i\in\ints{p}}{i\ge 2,\ \alpha_i=\alpha_{i-1}+1}\cup\{1\}
\]
(where, for technical reasons, $1$ is included).

The following observation follows readily from the definitions.

\begin{figure}[tb]
  \centering
  \newcommand{\mystyle}[1]{\footnotesize {#1}}
  \psfrag{0}{\mystyle{$0$}}
  \psfrag{1}{\mystyle{$1$}}
  \psfrag{*}{\mystyle{$\star$}}
  \subfloat[][]{\includegraphics[width =.15\textwidth]{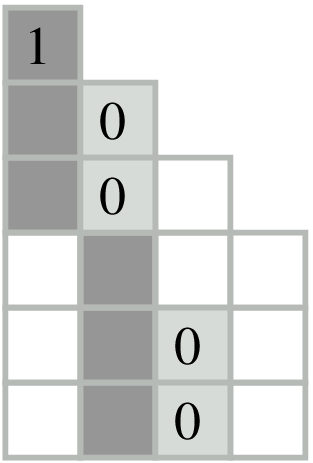}\label{sfl:zeroSC}}\hspace{10ex}
  \subfloat[][]{\includegraphics[width =.15\textwidth]{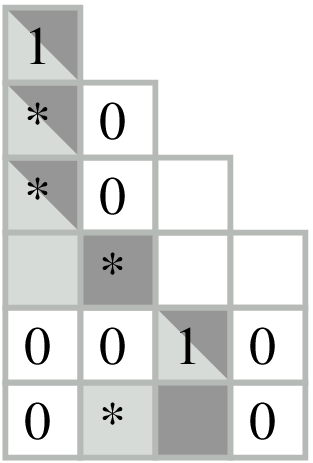}\label{sfl:feasiblePoint}}
  \caption[]{\subref{sfl:zeroSC}: Example for Remark~\ref{rem:SCIalpha}.
    Dark-gray entries indicate entries $(i,\alpha_i(I_0))$ and light-gray
    entries indicate the entries in $S_i(I_0)$ for $i = 6$.\\
    \subref{sfl:feasiblePoint}: Example for Lemma~\ref{lem:alpha:2}. As
    before, dark-gray entries indicate entries $(i,\alpha_i)$. Light-gray
    entries indicate entries $(i,\mu_i(I_0))$. The $\star$'s indicate $1$s set in the point
    $x^\star$ as constructed in Lemma~\ref{lem:alpha:2}.}
  \label{fig:ZeroSC_Feas}
\end{figure}

\begin{remark}\label{rem:SCIalpha}
  For each $i\in\ints{p}$ with $i\ge 2$ and $\alpha_{i}(I_0)<q(i)$, the set
  $S_i(I_0) \define \setdef{(k,\alpha_k(I_0)+1)}{k\in\ints{i}\setminus
    \Gamma(I_0)}$ is a shifted column with $S_i(I_0)\subseteq I_0$.
\end{remark}

Figure~\ref{fig:ZeroSC_Feas}~\subref{sfl:zeroSC} shows an example.

\begin{lemma}\label{lem:alpha:1}
  For each $i\in\ints{p}$, no vertex of $\orbipart{p}{q}\cap F$ has its
  $1$-entry in row~$i$ in a column $j\in \ints{q(i)}$ with
  $j>\alpha_i(I_0)$.
\end{lemma}

\begin{proof}
  Let $i \in \ints{p}$. We may assume $\alpha_i(I_0) < q(i)$, because
  otherwise the statement is trivially true.  Thus, $B \define
  \setdef{(i,j)\in\row{i}}{j>\alpha_i(I_0)} \neq \varnothing$.

  Let us first consider the case $i\in\Gamma(I_0)$. As we have
  $\alpha_i(I_0) < q(i) \leq i$ and $\alpha_1(I_0)=1$, there must be some
  $k<i$ such that $k \not \in \Gamma(I_0)$. Let~$k$ be maximal with this
  property. Thus, we have $k' \in \Gamma(I_0)$ for all $1<k< k'\le i$.
  According to Remark~\ref{rem:SCIalpha}, $x(B) - x(S_{k}(I_0)) \leq 0 $ is
  a shifted column inequality with $x(S_{k}(I_0))=0$, showing $x(B)=0$ as
  claimed in the lemma.

  Thus, let us suppose $i\in\ints{p}\setminus\Gamma(I_0)$. If
  $\alpha_i(I_0)\ge q(i)-1$, the claim holds trivially. Otherwise, $B' \define B
  \setminus \{(i,\alpha_i(I_0)+1)\} \neq \varnothing$. Similarly to the
  first case, now the shifted column inequality $x(B') - x(S_{i-1}(I_0))
  \leq 0$ proves the claim.
\end{proof}

For each $i \in \ints{p}$, we define $\mu_i(I_0) \define
\min\setdef{j\in\ints{q(i)}}{(i,j)\not\in I_0} $. Because of Property~P1,
the sets over which we take minima here are non-empty.

\begin{lemma}\label{lem:alpha:2}
  If we have $\mu_i(I_0) \leq \alpha_i(I_0)$ for all $i\in\ints{p}$, then
  the point $x^{\star} = x^{\star}(I_0) \in \{0,1\}^{\orbipartinds{p}{q}}$
  defined by $x^{\star}_{i,\alpha_i(I_0)} = 1$ for all $i \in \Gamma(I_0)$,
  $x^{\star}_{i,\mu_i(I_0)} = 1$ for all $i\in \ints{p}\setminus
  \Gamma(I_0)$, and all other components being zero, is contained in
  $\orbipart{p}{q}\cap F$.
\end{lemma}

\begin{proof}
  Due to $\alpha_i(I_0)\le\alpha_{i-1}(I_0)+1$ for all $i\in\ints{p}$ with
  $i\ge 2$, the point~$x^{\star}$ is contained in $\orbipart{p}{q}$.  It
  follows from the definitions that $x^{\star}$ does not have a $1$-entry
  at a position in~$I_0$. Thus, by Remark~\ref{rem:props}, we have
  $x^{\star}\in F$.
\end{proof}

\noindent Figure~\ref{fig:ZeroSC_Feas}~\subref{sfl:feasiblePoint} shows an
example for the point constructed in Lemma~\ref{lem:alpha:2}.

We now characterize the case $\orbipart{p}{q}\cap F=\varnothing$ (leading
to pruning the corresponding node in the branch-and-cut tree) and describe
the set~$I^{\star}_0$.

\begin{proposition}\label{prop:orbiFix}\
  \begin{enumerate}
  \item\label{prop:orbiFix:prune} We have $\orbipart{p}{q}\cap F=\varnothing$
    if and only if there exists $i\in\ints{p}$ with
    $\mu_i(I_0)>\alpha_i(I_0)$.
  \item\label{prop:orbiFix:notprune} If $\mu_i(I_0)\le \alpha_i(I_0)$ holds
    for all $i\in\ints{p}$, then the following is true.
    \begin{enumerate}
    \item\label{prop:orbiFix:notprune:fixone} For all $i\in\ints{p}\setminus
      \Gamma(I_0)$, we have
      \[
      I^{\star}_0\cap\row{i} = \{ (i,j) \in \row{i} \; : \: (i,j) \in I_0
      \text{ or }j > \alpha_i(I_0)\}.
      \]
    \item\label{prop:orbiFix:notprune:fixtwo} For all $i \in \ints{p}$ with
      $\mu_i(I_0) = \alpha_i(I_0)$, we have
      \[
      I^{\star}_0\cap\row{i} = \row{i}\setminus\{(i,\alpha_i(I_0))\}.
      \]
    \item\label{prop:orbiFix:notprune:fixthree}
      For all $s\in \Gamma(I_0)$ with $\mu_s(I_0)<\alpha_s(I_0)$ the
      following holds: If there is some $i \geq s$ with
      $\mu_i(I_0)>\alpha_i(I_0\cup\{(s,\alpha_s(I_0))\}) $, then we have
      \[
      I^{\star}_0\cap\row{s}\ =\ \row{s}\setminus\{(s,\alpha_s(I_0))\}.
      \]
      Otherwise, we have
      \[
      I^{\star}_0\cap\row{s}\ =\
      \setdef{(s,j)\in\row{s}}{(s,j)\in I_0\text{ or }j>\alpha_s(I_0)}.
      \]
    \end{enumerate}
  \end{enumerate}
\end{proposition}

\begin{figure}[tb]
  \centering
  \newcommand{\mystyle}[1]{\footnotesize {#1}}
  \psfrag{0}{\mystyle{$0$}}
  \psfrag{1}{\mystyle{$1$}}
  \subfloat[][]{\includegraphics[width =.15\textwidth]{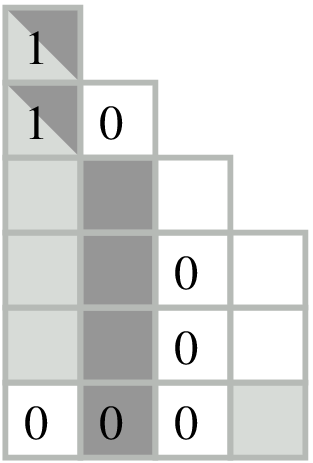}\label{sfl:infeasible1}}\hfill
  \subfloat[][]{\includegraphics[width =.15\textwidth]{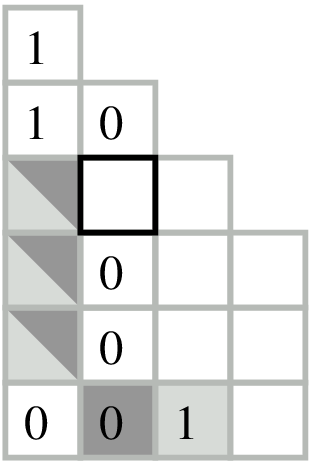}\label{sfl:fixOne}}\hfill
  \subfloat[][]{\includegraphics[width =.15\textwidth]{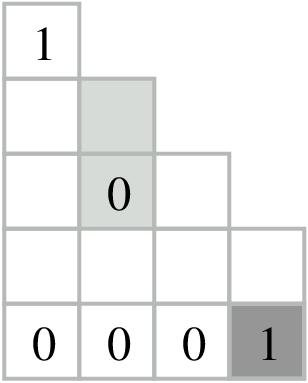}\label{sfl:seqfixSCICI}}\hfill
  \subfloat[][]{\includegraphics[width =.15\textwidth]{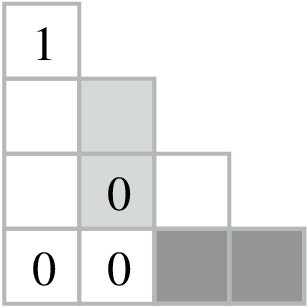}\label{sfl:seqfixOrbiFixSCI}}
  \caption[]{\subref{sfl:infeasible1}: Example for
    Proposition~\ref{prop:orbiFix}~\eqref{prop:orbiFix:prune}. Light-gray
    entries indicate the entries $(i,\mu_i(I_0))$ and dark-gray entries
    indicate entries $(i,\alpha_i(I_0))$. \subref{sfl:fixOne}: Example of
    fixing an entry to~$1$ for Proposition~\ref{prop:orbiFix}~(2c). As
    before light-gray entries indicate entries $(i,\mu_i(I_0))$. Dark-gray
    entries indicate
    entries $(i,\alpha_i(I_0\cup\{(s,\alpha_s(I_0))\}))$ with $s = 3$.\\
    \subref{sfl:seqfixSCICI} and \subref{sfl:seqfixOrbiFixSCI}: Gray
    entries show the SCIs used in the proofs of Parts~1(a) and~1(b) of
    Theorem~\ref{thm:seqfixOrbi}, respectively.}
  \label{fig:sympartopt}
\end{figure}

\begin{proof}
  Part~\ref{prop:orbiFix:prune} follows from Lemmas~\ref{lem:alpha:1}
  and~\ref{lem:alpha:2} (see also
  Figure~\ref{fig:sympartopt}~\subref{sfl:infeasible1}).

  In order to prove Part~\ref{prop:orbiFix:notprune}, let us assume that
  $\mu_i(I_0)\le \alpha_i(I_0)$ holds for all $i\in\ints{p}$.  For
  Part~\ref{prop:orbiFix:notprune:fixone}, let
  $i\in\ints{p}\setminus\Gamma(I_0)$ and $(i,j)\in\row{i}$. Due to
  $I_0\subseteq I^{\star}_0$, we only have to consider the case $(i,j)\not\in
  I_0$.  If $j>\alpha_i(I_0)$, then, by Lemma~\ref{lem:alpha:1}, we find
  $(i,j)\in I^{\star}_0$. Otherwise, the point that is obtained from
  $x^{\star}(I_0)$ (see Lemma~\ref{lem:alpha:2}) by moving the $1$-entry in
  position $(i,\mu_i(I_0))$ to position $(i,j)$ is contained in
  $\orbipart{p}{q}\cap F$, proving $(i,j)\not\in I^{\star}_0$.

  In the situation of Part~\ref{prop:orbiFix:notprune:fixtwo}, the claim
  follows from Lemma~\ref{lem:alpha:1} and $\orbipart{p}{q}\cap
  F\not=\varnothing$ (due to Part~1).

  For Part~\ref{prop:orbiFix:notprune:fixthree}, let $s\in\Gamma(I_0)$ with
  $\mu_s(I_0)<\alpha_s(I_0)$ and define $I_0' \define
  I_0\cup\{(s,\alpha_s(I_0))\}$. It follows that we have
  $\mu_i(I'_0)=\mu_i(I_0)$ for all $i\in\ints{p}$; compare also
  Figure~\ref{fig:sympartopt}~\subref{sfl:fixOne}.

  Let us first consider the case that there is some $i\ge s$ with
  $\mu_i(I_0)>\alpha_i(I_0') $. Part~\ref{prop:orbiFix:prune} (applied to
  $I'_0$ instead of~$I_0$) implies that $\orbipart{p}{q}\cap F$ does not
  contain a vertex~$x$ with $x_{s,\alpha_s(I_0)}=0$. Therefore, we have
  $(s,\alpha_s(I_0))\in I^{\star}_1$, and thus
  $I^{\star}_0\cap\row{s}=\row{s}\setminus\{(s,\alpha_s(I_0))\}$ holds
  (where for ``$\subseteq$'' we exploit $\orbipart{p}{q}\cap
  F\not=\varnothing$ by Part~\ref{prop:orbiFix:prune}, this time applied
  to~$I_0$).

  The other case of Part~\ref{prop:orbiFix:notprune:fixthree} follows from
  $s\not\in\Gamma(I_0')$ and $\alpha_s(I_0')=\alpha_s(I_0)-1$. Thus,
  Part~\ref{prop:orbiFix:notprune:fixone} applied to $I'_0$ and $s$ instead
  of $I_0$ and $i$, respectively, yields the claim (because of
  $(s,\alpha_s(I_0))\not\in I^{\star}_0$ due to $s\in\Gamma(I_0)$ and
  $\orbipart{p}{a}\cap F\not=\varnothing$).
\end{proof}

\subsection{Sequential Fixing for Partitioning Orbitopes}
\label{subsec:fixSCI}

Let us, for some fixed $p \geq q \geq 2$, denote by $\syssci$ the system of
the non-negativity inequalities, the row-sum equations (each one written as
two inequalities, in order to be formally correct) and all shifted column
inequalities. Thus, according to Theorem~\ref{thm:OrbitpartDesc},
$\orbipart{p}{q}$ is the set of all $x\in \R^{\orbipartinds{p}{q}}$ that
satisfy~$\syssci$.  Let $\sysci$ be the subsystem of $\syssci$ containing
only the column inequalities (and all non-negativity inequalities and
row-sum equations).

At first sight, it is not clear whether sequential fixing with the
exponentially large system $\syssci$ can be done efficiently. A closer look
at the problem reveals, however, that one can utilize the linear time
separation algorithm for shifted column inequalities (mentioned in
Sect.~\ref{sec:orbitopes}) in order to devise an algorithm for this
sequential fixing, whose running time is bounded by $\bigo{\varrho pq}$,
where $\varrho$ is the number of variables that are fixed by the procedure.

In fact, one can achieve more: One can compute sequential fixings with
respect to the affine hull of the partitioning orbitope. In order to
explain this, consider a polytope $P=\setdef{x\in\cube{d}}{Ax\le b}$, and
let $S\subseteq\R^d$ be some affine subspace containing~$P$. As before, we
denote the knapsack relaxations of~$P$ obtained from $Ax\le b$ by $P_1$,
\dots, $P_m$. Let us define $\fixfaceaff{P_r}{F}{S}$ as the smallest cube
face that contains $P_r\cap S\cap\{0,1\}^d\cap F$. Similarly to the
definition of $\fixface{Ax\le b}{F}$, denote by~$\fixfaceaff{Ax\le
  b}{F}{S}$ the face of~$\cube{d}$ that is obtained by setting $G \define
F$ and then iteratively replacing $G$ by $\fixfaceaff{P_r}{G}{S}$ as long
as there is some~$r\in\ints{m}$ with $\fixfaceaff{P_r}{G}{S}\subsetneq G$.
We call $\fixfaceaff{Ax\le b}{F}{S}$ the \emph{sequential fixing of $Ax\le
  b$ at~$F$ relative to~$S$}. Obviously, we have
$\fixface{P}{F}\subseteq\fixfaceaff{Ax\le b}{F}{S}\subseteq\fixface{Ax\le
  b}{F}$. In contrast to sequential fixing, sequential fixing relative to
affine subspaces \emph{in general} is $\NP$-hard (as it can be used to decide
whether a linear equation has a 0/1-solution).

\begin{theorem}\label{thm:seqfixOrbi}\
\begin{enumerate}
\item   There are cube faces $F^1$, $F^2$, $F^3$ with the following properties:
\begin{enumerate}
\item $\fixface{\syssci}{F^1} \subsetneq \fixface{\sysci}{F^1}$\label{part:SCI_CI}
\item $\fixfaceaff{\sysci}{F^2}{\aff{\orbipart{p}{q}}} \subsetneq
       \fixface{\syssci}{F^2}$\label{part:affCI_SCI}
\item $\fixfaceaff{\syssci}{F^3}{\aff{\orbipart{p}{q}}} \subsetneq
        \fixfaceaff{\sysci}{F^3}{\aff{\orbipart{p}{q}}}$\label{part:affSCI_affCI}
\end{enumerate}
\item For all cube faces~$F$, we have
  $\fixfaceaff{\syssci}{F}{\aff{\orbipart{p}{q}}}=\fixface{\orbipart{p}{q}}{F}$.\label{part:OF_SCI}
\end{enumerate}
\end{theorem}

\begin{proof}
  For Part~\eqref{part:SCI_CI}, we chose $p=5$, $q=4$, and define the cube
  face~$F^1$ via $I^1_0=\{(3,2),(5,1),(5,2),(5,3)\}$ and
  $I^1_1=\{(1,1),(5,4)\}$. The shifted column inequality with shifted
  column $\{(2,2),(3,2)\}$ and bar $\{(5,4)\}$ allows to fix $x_{22}$ to
  $1$ (see Fig.~\ref{fig:sympartopt}~\subref{sfl:seqfixSCICI}), while no
  column inequality (and no non-negativity constraint and no row-sum
  equation) allows to fix any variable.

  For Part~\eqref{part:affCI_SCI}, let $p=4$, $q=4$, and define~$F^2$ via
  $I^2_0=\{(3,2), (4,1), (4,2)\}$ and $I^2_1=\{(1,1)\}$. Exploiting that
  $x_{43}+x_{44}=1$ for all $x\in \aff{\orbipart{p}{q}}\cap F^2$, we can
  use the column inequality with column $\{(2,2),(3,2)\}$ and bar
  $\{(4,3),(4,4)\}$ to fix $x_{22}$ to one (see
  Fig.~\ref{fig:sympartopt}~\subref{sfl:seqfixOrbiFixSCI}), while no fixing
  is possible with $\syssci$ only.

  For Part~\eqref{part:affSCI_affCI}, we can use $F^3=F^1$.

  In order to prove Part~\eqref{part:OF_SCI}, we have to show for every cube
  face~$F$ that
  \[
  \fixfaceaff{\syssci}{F}{\aff{\orbipart{p}{q}}} \subseteq
  \fixface{\orbipart{p}{q}}{F}
  \]
  holds.  We use the notation introduced in Sect.~\ref{sec:orbifix}.  The
  crucial fact is that every point $x \in F \cap \aff{\orbipart{p}{q}}$
  satisfies $x(B)=1$ for every $B\subseteq\row{i}$ such that
  $\{(i,\mu_i(I_0)),\dots,(i,q(i))\}\subseteq B$.

  Let us first consider the case $\smash{\fixface{\orbipart{p}{q}}{F}} =
  \varnothing$. Due to Part~\ref{prop:orbiFix:prune} of
  Prop.~\ref{prop:orbiFix} there is some $i\in\ints{p}$ with
  $\mu_i(I_0)>\alpha_i(I_0)$. Therefore, the SCI
  $x(\tilde{B})-x(\tilde{S})\le 0$ constructed in the proof of
  Lemma~\ref{lem:alpha:1} with $(i,j)=(i,\mu_i(I_0))$ has $x(\tilde{B})=1$
  for all $x\in F\cap\aff{\orbipart{p}{q}}$, but $x(\tilde{S})=0$ due to
  $\tilde{S}\subseteq I_0$. This shows that we indeed have
  \[
  \fixfaceaff{\syssci}{F}{\aff{\orbipart{p}{q}}} = \varnothing
  \]
  in this case.

  Otherwise (i.e., $\fixface{\orbipart{p}{q}}{F} \neq \varnothing$), it
  suffices to show for each $(k,\ell) \in I^{\star}_0 \setminus I_0$ that
  there is some SCI that can only be satisfied by some $x \in F \cap
  \aff{\orbipart{p}{q}}$ if $x_{k\ell} = 0$ holds. Due to
  Part~\ref{prop:orbiFix:notprune} of Prop.~\ref{prop:orbiFix}, we have to
  consider two cases.\smallskip

  \noindent
  \emph{Case 1:} We have $\ell > \alpha_k(I_0)$. Then the SCI $x(\tilde{B})
  - x(\tilde{S}) \leq 0$ constructed in the proof of
  Lemma~\ref{lem:alpha:1} with $(i,j) = (k,\ell)$ implies $x(\tilde{B}) =
  0$ for all $x \in F$ (because of $\tilde{S} \subseteq I_0$), which yields
  $x_{k\ell} = 0$ due to $(k,\ell) \in \tilde{B}$.\smallskip

  \noindent
  \emph{Case 2:} We have $k \in \Gamma(I_0)$ with $\ell < \alpha_k(I_0)$
  and there is some $r \geq k$ with $\mu_r(I_0) > \alpha_r(I_0 \cup
  \{(k,\alpha_k(I_0))\})$. Then the SCI $x(\tilde{B}) - x(\tilde{S}) \leq
  0$ constructed in the proof of Lemma~\ref{lem:alpha:1} with $(i,j) =
  (r,\mu_r(I_0))$ (and $I_0$ replaced by $I_0 \cup \{(k, \alpha_k(I_0))\}$)
  satisfies, for each $x \in F \cap \aff{\orbipart{p}{q}}$,
  $x(\tilde{B})=1$ and $x(\tilde{S}) = x_{(k,\alpha_k(I_0))}$ (due to
  $\tilde{S} \subseteq I_0 \cup \{(k, \alpha_k(I_0))\}$), which implies
  $x_{(k,\alpha_k(I_0))} = 1$, and hence (as $x \in F \cap
  \aff{\orbipart{p}{q}}$) $x_{k\ell} = 0$, because of $\ell \neq
  \alpha_k(I_0)$.
\end{proof}

The different versions of sequential fixing for partitioning orbitopes are
dominated by each other in the following sequence:
\[
\sysci \rightarrow \{\syssci, \text{affine }\sysci\} \rightarrow \text{affine }\syssci,
\]
which finally is as strong as orbitopal fixing. For each of the arrows
there exists an instance for which dominance is strict.  The examples in
the proof of Theorem~\ref{thm:seqfixOrbi} also show that there is no
general relation between $\syssci$ and affine~$\sysci$.

In particular, we could compute orbitopal fixings by the polynomial time
algorithm for sequential fixing relative to $\aff{\orbipart{p}{q}}$. It
turns out, however, that this is not the preferable choice. In fact, we
will describe below a linear time algorithm for solving the orbitopal
fixing problem directly.

\subsection{An Algorithm for Orbitopal Fixing}

Algorithm~\ref{alg:OrbitopalFixing} describes a method to compute the
simultaneous fixing $(I^{\star}_0,I^{\star}_1)$ from $(I_0,I_1)$ (which are
assumed to satisfy Properties~P1 and~P2). Note that we use~$\beta_i$ for
$\alpha_i(I_0\cup\{(s,\alpha_s(I_0))\})$.

\begin{algorithm}[t]
  \floatname{algorithm}{\footnotesize Algorithm}
  \footnotesize
  \caption{\footnotesize Orbitopal Fixing}
  \label{alg:OrbitopalFixing}
  \begin{algorithmic}[1]
    \STATE Set $I^{\star}_0 \leftarrow I_0$, $I^{\star}_1 \leftarrow I_1$, $\mu_1\leftarrow 1$,
    $\alpha_1 \leftarrow 1$, and $\Gamma = \varnothing$.
    \FOR{$i = 2, \dots, p$}\label{step:ForAlpha}
    \STATE compute $\mu_i \leftarrow \min \setdef{j}{(i,j) \not\in I_0}$.\label{step:mui}\\
    \IF{$\alpha_{i-1} = q(i)$ or $(i,\alpha_{i-1}+1) \in I_0$}
    \STATE $\alpha_i \leftarrow \alpha_{i-1}$
    \ELSE
    \STATE $\alpha_i \leftarrow \alpha_{i-1}+1$, $\Gamma \leftarrow \Gamma \cup \{i\}$
    \ENDIF
    \IF{$\mu_i>\alpha_i$}
      \STATE return ``Orbitopal fixing is empty''\label{step:infeasible}
    \ENDIF
    \STATE Set $I^{\star}_0 \leftarrow I^{\star}_0 \cup
    \setdef{(i,j)}{j>\alpha_i}$.\label{step:ZeroFixing}
    \IF{$|I^{\star}_0\cap\row{i}|=q(i)-1$} \label{step:Istarone}
    \STATE set $I^{\star}_1 \leftarrow I^{\star}_1\cup(\row{i}\setminus I^{\star}_0)$.
    \label{step:ForAlphaEnd}
    \ENDIF
    \ENDFOR
    \FORALL{$s\in \Gamma$ with $(s, \alpha_s) \notin I^{\star}_1$}\label{step:SLoop}
    \STATE Set $\beta_s \leftarrow \alpha_s - 1$.\\
    \FOR{$i = s+1, \dots, p$}\label{step:ForBeta}
    \IF{$\beta_{i-1} = q(i)$ or $(i, \beta_{i-1}+1) \in I_0$}
    \STATE $\beta_i \leftarrow \beta_{i-1}$
    \ELSE
    \STATE $\beta_i \leftarrow \beta_{i-1} + 1$
    \ENDIF
    \IF{$\mu_i > \beta_i$}\label{step:ifmubeta}
    \STATE $I^{\star}_1 \leftarrow I^{\star}_1 \cup \{(s,\alpha_s)\}$ and
    $I^{\star}_0 \leftarrow \row{s} \setminus \{(s,\alpha_s)\}$.\label{step:OneFixing}\\
    \STATE Proceed with the next $s$ in Step~\ref{step:SLoop}.
    \label{step:ForBetaEnd}
    \ENDIF
    \ENDFOR
    \ENDFOR
  \end{algorithmic}
\end{algorithm}

\begin{theorem}\label{thm:algorbifix}
  The orbitopal fixing problem for partitioning orbitopes can be solved in
  time $\bigo{pq}$ (by a slight modification of
  Algorithm~\ref{alg:OrbitopalFixing}).
\end{theorem}
\begin{proof}
  The correctness of the algorithm follows from the structural results
  given in Proposition~\ref{prop:orbiFix}.

  In order to prove the statement on the running time, let us assume that
  the data structures for the sets~$I_0$, $I_1$, $I^{\star}_0$,
  and~$I^{\star}_1$ allow both membership testing and addition of single
  elements in constant time (e.g., the sets can be stored as bit vectors).

  As none of the Steps~\ref{step:mui} to~\ref{step:ForAlphaEnd} needs more
  time than $\bigo{q}$, we only have to take care of the second part of the
  algorithm starting in Step~\ref{step:SLoop}. (In fact, used verbatim as
  described above, the algorithm might need time $\Omega(p^2)$.)

  For $s,s'\in \Gamma$ with $s<s'$ denote the corresponding $\beta$-values
  by $\beta_i$ ($i\ge s$) and by~$\beta'_i$ ($i\ge s'$), respectively. We
  have $\beta_i\le\beta'_i$ for all $i\ge s'$, and furthermore, if equality
  holds for one of these~$i$, we can deduce $\beta_k=\beta'_k$ for all
  $k\ge i$. Thus, as soon as a pair $(i,\beta_i)$ is used a second time in
  Step~\ref{step:ifmubeta}, we can break the for-loop in
  Step~\ref{step:ForBeta} and reuse the information that we have obtained
  earlier.

  This can, for instance, be organized by introducing, for each
  $(i,j)\in\orbipartinds{p}{q}$, a flag $f(i,j)\in\{\red,\green,\white\}$
  (initialized by $\white$), where $f(i,j)=\red/\green$ means that we have
  already detected that $\beta_i=j$ eventually leads to a positive/negative
  test in Step~\ref{step:ifmubeta}.  The modifications that have to be
  applied to the second part of the algorithm are the following: The
  selection of the elements in~$\Gamma$ in Step~\ref{step:SLoop} must be
  done in increasing order.  Before performing the test in
  Step~\ref{step:ifmubeta}, we have to check whether $f(i,\beta_i)$ is
  $\green$. If this is true, then we can proceed with the next~$s$ in
  Step~\ref{step:SLoop}, after setting all flags $f(k,\beta_k)$ to $\green$
  for $s \leq k < i$. Similarly, we set all flags $f(k,\beta_k)$ to $\red$
  for $s \leq k \leq i$, before switching to the next~$s$ in
  Step~\ref{step:ForBetaEnd}. And finally, we set all flags $f(k,\beta_k)$
  to $\green$ for $s\le k\le p$ at the end of the body of the $s$-loop
  starting in Step~\ref{step:SLoop}.

  As the running time of this part of the algorithm is proportional to the
  number of flags changed from $\white$ to $\red$ or $\green$, the total
  running time indeed is bounded by $\bigo{pq}$ (since a flag is never
  reset).
\end{proof}

\section{Fixing for Packing and Covering Orbitopes}
\label{sec:packcov}

The packing orbitope~$\orbipack{p}{q}$ obviously can be obtained from the
partitioning orbitope~$\orbipart{p+1}{q+1}$ by projecting out the first
column and row, i.e., by orthogonal projection to the coordinate subspace
associated with
\[
\setdef{(i,j) \in \orbipartinds{p+1}{q+1}}{i,j > 1}
\]
(and renaming the variables appropriately), see also~\cite{KaPf08}.

In general, for $J \subseteq \ints{d}$, the orthogonal projection
$\pi:\R^{\ints{d}} \rightarrow \R^J$, a polytope $P \subseteq \cube{d}$,
and some face~$F$ of the cube $\cube{J} \define [0,1]^J$, we have
\[
\fixface{F}{\pi(P)}=\pi(\fixface{\pi^{-1}(F)}{P}),
\]
since for every face $G$ of~$\cube{J}$
\[
(\pi(P)\cap F)\subseteq G\Longleftrightarrow (P\cap\pi^{-1}(F))\subseteq\pi^{-1}(G)
\]
holds (simply because taking preimages commutes with taking intersections).

Thus, the following result for packing orbitopes follows readily from
Theorem~\ref{thm:algorbifix}.

\begin{corollary}\label{cor:PackingOrbitopes}
  Variable fixing for packing orbitopes~$\orbipack{p}{q}$ can be done in
  time $\bigo{pq}$ by reduction to orbitopal fixing for
  $\orbipart{p+1}{q+1}$.
\end{corollary}

In contrast to this, variable fixing for covering
orbitopes~$\orbicov{p}{q}$ cannot be done in polynomial time, unless
$\polytime = \NP$, as the following result implies.

\begin{figure}[tb]
  \begin{minipage}{0.4\textwidth}
    \newcommand{\mystyle}[1]{\footnotesize {#1}}
    \psfrag{1}{\mystyle{$2$}}
    \psfrag{2}{\mystyle{$4$}}
    \psfrag{3}{\mystyle{$6$}}
    \psfrag{4}{\mystyle{$8$}}
    \psfrag{5}{\mystyle{$10$}}
    \psfrag{6}{\mystyle{$12$}}
    \psfrag{7}{\mystyle{$14$}}
    \psfrag{8}{\mystyle{$16$}}
    \psfrag{9}{\mystyle{$18$}}
    \psfrag{10}{\mystyle{$20$}}
    \psfrag{11}{\mystyle{$22$}}
    \psfrag{12}{\mystyle{$24$}}
    \psfrag{13}{\mystyle{$26$}}
    \psfrag{14}{\mystyle{$28$}}
    \includegraphics[width=\textwidth]{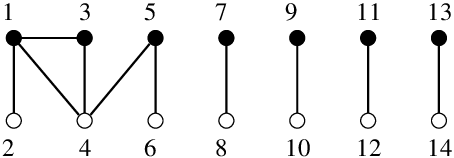}
  \end{minipage}
  \hfill
  \begin{minipage}{0.45\textwidth}
    \newcommand{\mystyle}[1]{\footnotesize {#1}}
    \psfrag{a}{\mystyle{$4$}}
    \psfrag{b}{\mystyle{$2$}}
    \psfrag{c}{\mystyle{$6$}}
    \psfrag{d}{\mystyle{$1$}}
    \psfrag{e}{\mystyle{$3$}}
    \psfrag{f}{\mystyle{$5$}}
    \psfrag{g}{\mystyle{$7$}}
    \includegraphics[width=\textwidth]{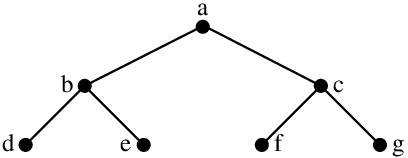}
  \end{minipage}
  \vspace{1ex}

  \begingroup
  \centering
  \tiny
  \newcolumntype{L}{@{\;\;}r}
  \newcommand{\1}{\textbf{1}}
  \begin{tabular}{c|rLLLLLLLLLLLLLLLLLLLLLLLLLLL|}\cline{2-29}
    & \1 & \1 & \1 & \1 & \1 & \1 & \1 & \1 & \1 & \1 & \1 & \1 & \1 &   &   &   &   &   &   &   &   &   &   &   &   &   &   &   \\
    & \1 & \1 & \1 & \1 & \1 &   &   &   &   &   &   &   &   & \1 & \1 & \1 & \1 & \1 & \1 & \1 & \1 &   &   &   &   &   &   &   \\
    & \1 &   &   &   &   & \1 & \1 & \1 & \1 &   &   &   &   & \1 & \1 & \1 & \1 &   &   &   &   & \1 & \1 & \1 & \1 &   &   &   \\
    \cline{2-29}
      \{2,4\} & 0 &\1 & 0 &   & 0 & 0 & 0 & 0 & 0 & 0 & 0 & 0 & 0 & 0 & 0 & 0 & 0 & 0 & 0 & 0 & 0 & 0 & 0 & 0 & 0 & 0 & 0 & 0 \\
      \{2,6\} & 0 &\1 & 0 & 0 & 0 &   & 0 & 0 & 0 & 0 & 0 & 0 & 0 & 0 & 0 & 0 & 0 & 0 & 0 & 0 & 0 & 0 & 0 & 0 & 0 & 0 & 0 & 0 \\
      \{2,8\} & 0 &\1 & 0 & 0 & 0 & 0 & 0 &   & 0 & 0 & 0 & 0 & 0 & 0 & 0 & 0 & 0 & 0 & 0 & 0 & 0 & 0 & 0 & 0 & 0 & 0 & 0 & 0 \\
      \{6,8\} & 0 & 0 & 0 & 0 & 0 &\1 & 0 &   & 0 & 0 & 0 & 0 & 0 & 0 & 0 & 0 & 0 & 0 & 0 & 0 & 0 & 0 & 0 & 0 & 0 & 0 & 0 & 0 \\
     \{8,10\} & 0 & 0 & 0 & 0 & 0 & 0 & 0 &   & 0 &\1 & 0 & 0 & 0 & 0 & 0 & 0 & 0 & 0 & 0 & 0 & 0 & 0 & 0 & 0 & 0 & 0 & 0 & 0 \\
    \{10,12\} & 0 & 0 & 0 & 0 & 0 & 0 & 0 & 0 & 0 &\1 & 0 &   & 0 & 0 & 0 & 0 & 0 & 0 & 0 & 0 & 0 & 0 & 0 & 0 & 0 & 0 & 0 & 0 \\
    \{14,16\} & 0 & 0 & 0 & 0 & 0 & 0 & 0 & 0 & 0 & 0 & 0 & 0 & 0 &\1 & 0 &   & 0 & 0 & 0 & 0 & 0 & 0 & 0 & 0 & 0 & 0 & 0 & 0 \\
    \{18,20\} & 0 & 0 & 0 & 0 & 0 & 0 & 0 & 0 & 0 & 0 & 0 & 0 & 0 & 0 & 0 & 0 & 0 &\1 & 0 &   & 0 & 0 & 0 & 0 & 0 & 0 & 0 & 0 \\
    \{22,24\} & 0 & 0 & 0 & 0 & 0 & 0 & 0 & 0 & 0 & 0 & 0 & 0 & 0 & 0 & 0 & 0 & 0 & 0 & 0 & 0 & 0 &\1 & 0 &   & 0 & 0 & 0 & 0 \\
    \{26,28\} & 0 & 0 & 0 & 0 & 0 & 0 & 0 & 0 & 0 & 0 & 0 & 0 & 0 & 0 & 0 & 0 & 0 & 0 & 0 & 0 & 0 & 0 & 0 & 0 & 0 &\1 & 0 &   \\
    \cline{2-29}
    \multicolumn{1}{c}{} & \multicolumn{1}{L}{\hphantom{1}1} & \hphantom{1}2 & \hphantom{1}3 & \hphantom{1}4 & \hphantom{1}5 & \hphantom{1}6 &
    \hphantom{1}7 & \hphantom{1}8 & \hphantom{1}9 & 10 & 11 & 12 & 13 & 14 & 15 & 16 & 17 & 18 & 19 & 20 &
    21 & 22 & 23 & 24 & 25 & 26 & 27 & \multicolumn{1}{@{}r}{28}\\
  \end{tabular}
  \endgroup
  \caption{Example for the construction in the proof of
    Theorem~\ref{thm:orbicov}. Let $G$ be the depicted graph on the top left,
    and let $k = 7$, i.e., $\kappa = \lceil\log_2(k+1)\rceil
    = 3$ and $\tilde{k} = 7$. Thus, no new edge is needed, and we have $G =
    \tilde{G}$. In the matrix, $0$s correspond to elements of $I_0$, $\textbf{1}$s
    have been set in the construction of a feasible solution, and empty
    entries correspond to $0$s set in the construction. The top right shows
    the binary tree used for the construction with $a_1, \dots, a_7 =
    2, 6, 10, 14, 18, 22, 26$.}
  \label{fig:CoveringExample}
\end{figure}

\begin{theorem}\label{thm:orbicov}
  The problem to decide whether, for given $I_0 \subseteq \ints{p} \times
  \ints{q}$, the covering orbitope~$\orbicov{p}{q}$ contains a
  vertex~$x^\star \in \orbicov{p}{q}$ with $x^\star_{ij} = 0$ for all
  $(i,j)\in I_0$, is $\NP$-complete.
\end{theorem}

\begin{proof}
  It suffices to show that one can construct, for each graph $G = (V,E)$
  and $k \in \N$ with $k \leq \card{V}$, (in time bounded polynomially
  in~$\card{V}$) an instance of the decision problem described in the
  theorem whose answer is ``yes'' if and only if~$G$ has a vertex cover of
  size at most~$k$.

  Towards this end, let~$\kappa \define \lceil\log_2(k+1)\rceil$ be the
  smallest integer such that we have $\tilde{k} \define 2^{\kappa} - 1 \geq
  k$. Construct a graph~$\tilde{G} = (\tilde{V}, \tilde{E})$ by adding
  $\tilde{k}-k$ new edges (forming a matching) on $2(\tilde{k} - k)$ nodes
  disjoint from~$V$.  Thus~$\tilde{G}$ is a graph with $\card{\tilde{V}} =
  \card{V} + 2(\tilde{k} - k)$ nodes and $m \define \card{\tilde{E}} =
  \card{E} + \tilde{k} - k$ edges that has a vertex cover of size at
  most~$\tilde{k}$ if and only if~$G$ has a vertex cover of size at
  most~$k$.

  For the instance of the decision problem described in the theorem, let $p
  \define \kappa + m$, $q \define 2 \card{\tilde{V}}$, and assume
  $\tilde{V} = \{2,4,6, \dots, q\}$. Numbering the edges of~$\tilde{G}$ by
  \[
  \tilde{E} = \{e_1, e_2, \dots, e_m\} \text{ with } e_h = \{v_h,w_h\}
  \subseteq \tilde{V} \text{ for all } h \in \ints{m},
  \]
  we set
  \[
  I_0 \define \setdef{(\kappa + h, j)}{h \in \ints{m},\; j \in \ints{q}
    \setminus \{v_h, w_h\}}.
  \]
  See Figure~\ref{fig:CoveringExample} for an example.

  In order to prove that the answer to the constructed instance is ``yes''
  if and only if~$\tilde{G}$ has a vertex cover of size at
  most~$\tilde{k}$, let us call, for  $x^\star \in
  \{0,1\}^{\ints{p}\times\ints{q}}$, a pair~$(i,j)$ an \emph{alibi} (for
  column~$j$ of~$x^\star$), if $x_{i,j-1}^\star = 1$ and $x_{ij}^\star = 0$
  hold.

  If $x^\star \in \orbicov{p}{q}$ is a vertex of
  the covering orbitope~$\orbicov{p}{q}$ (i.e., a 0/1-point in the
  orbitope) with $x_{ij}^\star = 0$ for all $(i,j) \in I_0$, then
  \[
  C \define \setdef{v \in \tilde{V}}{x_{\kappa+h,v}^\star = 1 \text{ for
      some }h \in \ints{m}}
  \]
  is a vertex cover in~$\tilde{G}$ (due to $x^{\star}(\row{\kappa+h}) \geq
  1$ for all $h \in \ints{m}$). Moreover, for every $v \in C$, there is an
  alibi $(i,v)$ in some row~$i \in \ints{\kappa}$, since column~$v$
  of~$x^\star$ is lexicographically not larger than column~$v-1$. Again due
  to the lexicographical ordering of the columns, every vertex of the
  orbitope can have at most $2^{i-1}$ alibis in row~$i$. It follows that
  \[
  \card{C} \leq \sum_{i=1}^{\kappa}2^{i-1} = \sum_{i=0}^{\kappa-1}2^i = 2^{\kappa}-1 = \tilde{k}.
  \]

  Conversely, suppose $C \subseteq \tilde{V}$ is a vertex cover with
  $\card{C} \leq \tilde{k}$. We construct a 0/1-point $x^\star \in
  \orbicov{p}{q}$ with $x_{ij}^\star = 0$ for all $(i,j) \in I_0$ as
  follows. First, for each $h \in \ints{m}$, we set
  \[
  x_{\kappa+h,v_h}^\star \define
  \begin{cases}
    1 & \text{if } v_h \in C\\
    0 & \text{otherwise}
  \end{cases}
  \qquad
  \text{and}
  \qquad
  x_{\kappa+h,w_h}^\star \define
  \begin{cases}
    1 & \text{if } w_h \in C\\
    0 & \text{otherwise}.
  \end{cases}
  \]
  Since~$C$ is a vertex cover of~$\tilde{G}$, the part of~$x^\star$ that we
  have already constructed has at least one $1$-entry in every row
  $\kappa+1, \dots, p$. It thus remains to construct the first~$\kappa$
  rows such that they contain an alibi for every column~$v$ with~$v \in C$
  (and such that each of these rows contains at least one $1$-entry). This
  can, e.g., be achieved as follows.

  First, choose an arbitrary sequence $(a_1, \dots, a_{\tilde{k}})$ (of
  length $\tilde{k} = 2^{\kappa}-1$) of numbers in $\tilde{V} = \{2, 4, 6,
  \dots, q\}$ with $C \subseteq \{a_1, \dots, a_{\tilde{k}}\}$ (which is
  possible due to $\card{C} \leq \tilde{k}$).

  Then a complete rooted binary tree (embedded into the plane) of
  height~$\kappa$ (having $2^{\kappa}-1 = \tilde{k}$ nodes) is constructed
  in which the nodes receive pairwise different labels  $1, 2, \dots,
  \tilde{k}$. Furthermore, the labels have to be assigned in such a way
  that for every node labeled~$t$, we have $a_{\ell} \leq a_t$ for all
  labels $\ell$ in the \emph{left} subtree and $a_t \leq a_r$ for all
  labels $r$ in the \emph{right} subtree rooted at~$t$.

  Then we complete~$x^\star$ to a vertex of~$\orbicov{p}{q}$ by putting in
  each row $i \in \ints{\kappa}$ alibis at all positions $(i, a_t)$ for~$t$
  running through all labels of nodes at distance~$i-1$ from the root of
  the tree and filling the remaining components of~$x^\star$ accordingly.
\end{proof}

Of course, Theorem~\ref{thm:orbicov} implies that optimization over
covering orbitopes is $\NP$-hard. In particular, in contrast to the packing
and partitioning orbitopes, we cannot expect to find a tractable linear
description of~$\orbicov{p}{q}$, unless $\NP = \coNP$.

In fact, using ideas of the proof of Theorem~\ref{thm:orbicov}, one can
also establish other similar statements, in which $\orbicov{p}{q}$ is
replaced by the convex hull of all 0/1-matrices (whose columns are in
lexicographically non-increasing order) with \emph{at least}~$k$
one-entries per row for each $k\ge 1$, or with \emph{exactly} ~$k$
one-entries per row for each $k\ge 2$.

\section{Comparison with Isomorphism Pruning and Orbital Branching}
\label{sec:Margot}

In \cite{Mar02,Mar03b} Margot developed  a related, but more general approach to symmetry
breaking, called \emph{isomorphism pruning}. The main components are a setting rule for variables and a pruning
rule for nodes in a branch-and-bound tree to avoid consideration of
equivalent (partial) solutions. In this section we outline the differences
and similarities between Margot's and our approach when
specialized to the type of symmetries addressed by partitioning orbitopes.


Isomorphism pruning deals with arbitrary symmetries in any binary program
$\min\{c^T x : Ax \leq b, x \in \{0,1\}^d\} $, or even integer program
\cite{Mar03b}.  Let $G$ be a group of permutations of the variables
(inducing permutations of the components of~$c$, the columns of~$A$, and
the set of feasible solutions) such that for every $g \in G$ we have $g(c)
= c$ and $g(A) = \sigma_g(A)$, $g(b) = \sigma_g(b)$ for some permutation
$\sigma_g$ of the rows of $A$ resp.\ components of~$b$.  In particular, $G$
acts on the set of feasible solutions (via coordinate permutations) with
the property that the objective function $c^Tx$ is constant on every orbit.
Given an order of the variables (a \emph{rank vector}~$R$), Margot's
approach then assures that only (partial) solutions that are
lexicographically minimal in their orbit under $G$ are explored in the
branch-and-bound tree.  More precisely, a partial solution in some
branch-and-bound node~$N$ is identified with the sets $I^N_1$ and $I^N_0$
of variables that have been fixed to $1$ and~$0$, respectively, in the path
from the root to $N$. A partial solution is lex-min in its orbit, i.e., it
is a \emph{representative}, if the set $I^N_1$ is lexicographically minimal
with respect to the rank vector~$R$.  Thus, $N$ can be pruned if it is not
a representative.

Note that the definition of the lexicographical order relies upon a
particular total order of the variables defined by the rank vector~$R$. It
should be mentioned that in~\cite{Mar03b} this limitation was relaxed by
using an arbitrary order that can be determined during the branch-and-bound
process. Still, the same rank vector has to be used throughout.
Essentially, this means that whenever a branching is to be performed at
some level of the branch-and-bound tree for the first time, the branching
variable for this level can be freely chosen, and the rank vector is
extended.  In~\cite{Ost09} (see also~\cite{Mar10}), the restriction of a
global rank vector is dropped, too. The rank vector at some node $N$ of the
branch-and-bound tree is now given by the branching decisions from the root
node to $N$.

At every node $N$ of the branch-and-bound tree, two \emph{$0$-setting
  operations} are performed. If $N$ was created by fixing the variable
$x_f$ to $0$ then all variables in the (sub)orbit of $x_f$ under the
stabilizer of~$I^N_1$ are set to $0$, too, since a $1$ for any of these
variables would lead to a partial solution lex-greater than $I^N_1$.
Furthermore, in case the index~$h$ of the next branching variable is known
(e.g., if a global rank vector is available), it is repeatedly tested
whether any representative can be reached from~$N$ by checking whether the
current representative together with~$h$ is also a representative, i.e.,
whether $I^N_1 \cup \{h\}$ is lexicographically minimal under $G$.


For a comparison with orbitopal fixing, we consider isomorphism pruning
specialized for partitioning problems with 0/1-variables $x_{ij}$
satisfying $\sum_{j}x_{ij}=1$ for all~$i$, where the symmetry group $G$ is
assumed to be the group of all permutations of the columns of the variable
matrix~$x$.  Note that such problems do not require the elaborate machinery
of group theoretic algorithms developed by Margot for the general case of
more complicated symmetry groups.

We assume that the \emph{canonical rank vector~$R$} is used, i.e., the one
that describes the row-wise ordering of the variables $x_{ij}$.  For this
ordering, the representatives used by isomorphism pruning are in one-to-one
correspondence with faces of the cube $[0,1]^{\orbipartinds{p}{q}}$ having
nonempty intersections with $\orbipart{p}{q}$.

The row-wise ordering of the variables is natural choice. It moreover
turned out from our computational experiments with graph partitioning
problems that we could not find any alternative ordering yielding better
results for isomorphism pruning, at least for this application. This even
holds true for all variants we tested without a global rank vector (see
Sect.~\ref{sec:compu}).

\begin{figure}
	\centering
        \newcommand{\mystyle}[1]{\footnotesize {#1}}
        \psfrag{0}{\mystyle{$0$}}
        \psfrag{1}{\mystyle{$1$}}
	\includegraphics[height=3cm]{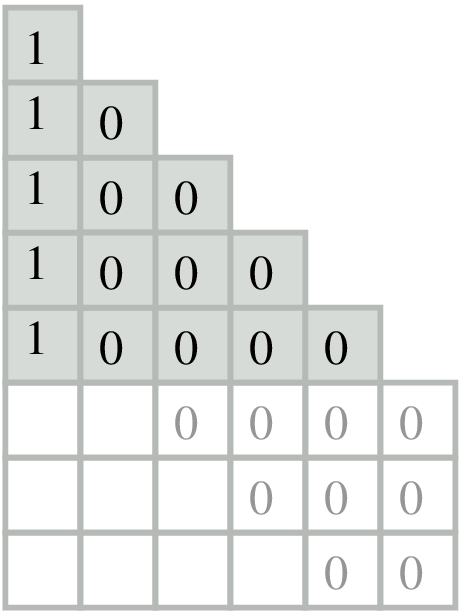}
	\caption{Example of a branch-and-bound node in which the variables in the first five rows have been fixed by minimum index branching. The zeroes below the fifth row have been set by orbitopal fixing after processing the current node, but not by isomorphism pruning.}
	\label{fig:fixmehr}
\end{figure}

If branch-and-bound trees for both methods are obtained by minimum index
branching, orbitopal fixing can be well compared to isomorphism pruning.
Indeed (provided that the nodes are also processed in the same order and no
cutting planes are added), orbitopal fixing will visit only
branch-and-bound nodes that isomorphism pruning visits as well. In any of
these nodes the orbitopal fixing algorithm does not perform $1$-fixings in
loop~\ref{step:SLoop} in Algorithm~\ref{alg:OrbitopalFixing}. However, it
may do more zero-fixings in loop~\ref{step:ForAlpha} than isomorphism
pruning. For an example see Figure~\ref{fig:fixmehr}.  Thus, the main
advantage of orbitopal fixing for our special case of symmetry can be seen
as deriving as early as possible conclusions that hold at every child node.

If branching rules different from minimum index branching are applied, the
trees produced by using orbitopal fixing and isomorphism pruning are not
comparable, in the sense that, in general, both variants will visit
branch-and-bound nodes (partial solutions) not visited by the other one.
But even at nodes that appear in both trees the behavior of the two
methods cannot really be compared unless both use the same representatives.
This, however, is only the case if the isomorphism pruning variant does
minimum index branching, yielding (locally) the canonical rank vector.
Hence, we are in the situation already discussed above (no matter by which
branching rule the orbitopal fixing variant has arrived at the node).


Another approach for avoiding symmetrical solutions is orbital
branching~\cite{LiOsRoSm06}.  It handles the same type of symmetry as
isomorphism pruning, but in a local manner: the symmetry group of the
current LP is computed on-the-fly for every node of the branch-and-bound
tree after removing fixed variables and inequalities that are satisfied
regardless of the unfixed variables.  The set of unfixed variables then
decomposes into orbits of equivalent variables under the symmetry group of
the current LP.  Then for some orbit of equivalent variables a two-way
branching is done, where in one branch the case is considered that all
variables in the orbit are zero, and in the other branch a chosen variable
from the orbit is fixed to one.  Orbital branching is thus not comparable
to either orbitopal fixing nor isomorphism pruning, as it can consider
symmetries that only arise after fixing variables. On the other hand a
problem usually looses global symmetry when variables are removed. As the
authors of~\cite{LiOsRoSm06} already point out, it is quite time consuming
to compute the symmetry group for every branch-and-bound node.
In~\cite{LiOsRoSm09} the authors report that a variant of orbital branching
exploiting global symmetry is computationally superior, i.e., instead of
computing the symmetry group of every local LP in the branch-and-bound
tree, the global symmetry group of the root LP is used. In this variant,
the symmetry group used in a branch-and-bound node is the subgroup of the
global symmetry group that setwise stabilizes the variables already fixed
to one. This, however, is very similar to isomorphism pruning.  If there is
some orbit $O$ of equivalent variables at some node $N$, then orbital
branching on $O$ is the same as pruning by isomorphism and $0$-setting
\emph{relative to the branching decisions leading to $N$}.  In this sense
orbital branching with global symmetry can be seen as isomorphism pruning
with local rank vectors.

In the case of partitioning orbitopes, the symmetry has very simple
structure.  Note that any symmetry considered here stems from the
permutation of columns of a matrix variable, which is a reasonable
restriction in particular with regard to the example application presented
in our graph partitioning formulation~\eqref{eq:kpart}.
At every branch-and-bound node, the columns of the matrix variable~$x$
decompose into one set of columns that are fixed elementwise, and one set
of columns which still can be permuted arbitrarily.

The special structure of partitioning problems considered here implies that
the orbits described above are the same in both variants of orbital
branching, at least as long as symmetry among rows of the matrix variable
(graph automorphisms in an instance of graph partitioning) are ignored.

The main advantage of orbital branching is its flexibility, e.g.,
orbital branching can be used at any node regardless how branching was
performed on other nodes of the branch-and-bound tree.

\section{Computational Experiments}
\label{sec:compu}

We performed computational experiments for the graph partitioning problem
mentioned in the introduction. The code is based on the SCIP~1.2.0
branch-and-cut framework~\cite{SCIP}, originally developed by
Achterberg~\cite{Ach07}. We use CPLEX 11.00 as the underlying LP solver.
The computations were performed on a 3.2 GHz Pentium~4 machine with 4~GB of
main memory and 2~MB cache running Linux. All computation times are CPU
seconds and are subject to a \emph{time limit of four hours}. Since in this
paper we are not interested in the performance of heuristics, we
initialized all computations with the \emph{optimal primal solution}.

We compare different variants of the code by counting \emph{winning}
instances. An instance is a winner for variant~A compared to variant~B,
if~A finished within the time limit and~B did not finish or needed a larger
CPU time; if~A did not finish, then the instance is a winner for~A in case
that~B did also not finish, leaving, however, a larger gap than~A. If the
difference between the times or gaps are below~1~sec. and 0.1~\%,
respectively, the instance is not counted.

In all variants, we fix the variables~$x_{ij}$ with~$j > i$ to zero.
Furthermore, we heuristically separate general clique inequalities
\[
\sum_{i,j \in C, i \neq j} y_{ij} \geq b,
\]
where
\[
b = \frac{1}{2}t(t-1)(q-r) + \frac{1}{2}t(t+1)r
\]
and $C \subseteq V$ is a clique of size~$tq+r > q$ with integers $t \geq
1$, $0\leq r < q$ (see \cite{ChRa93}).  The separation heuristic for a
fractional point~$y^{\star}$ follows ideas of Eisenbl\"atter~\cite{Eis01}.
We generate the graph $G' = (V, E')$ with $\{i,k\} \in E'$ if and only if
$\{i,k\} \in E$ and $y_{ik}^\star < b(b+1)/2$, where $y^\star$ is the
$y$-part of an LP solution. We search for maximum cliques in~$G'$ with the
specialized branch-and-bound method implemented in SCIP (with a node limit
of 10\,000), as well as with simple tabu search and greedy strategies. We
then check whether the corresponding inequality is violated.  We also
separate triangle inequalities and both kinds of cycle inequalities as
given in~\cite{ChRa93}.

After extensive testing, we decided to branch by default on the \emph{first
  index}, i.e., we branch on the first fractional $x$-variable in the
row-wise variable order used for defining orbitopes. A side-effect of this
choice is that this branching rule makes orbitopal fixing more comparable
to isomorphism pruning, in particular, to the variant using an \emph{a priori}
fixed variable order. It should be noted that this branching rule is
superior only when the vertices are ordered (i.e., the rows of~$x$ are
permuted) as follows: sort the vertices in descending order of their star
weight, i.e., the sum of the weights of incident edges.

We generated 36 random instances with~$p=40$ vertices and~$m$ edges of the
following types. We used $m=360$ (\emph{sparse}), $540$ (\emph{medium}),
and $720$ (\emph{dense}). For each type, we generated three instances by
picking edges uniformly at random (without recourse) until the specified
number of edges is reached. The edge weights are drawn independently and
uniformly at random from $[1000]$. For each instance we computed results
for $q =$ 3, 6, 9, and~12.

\begin{table}[tb]
  \caption{Results of the branch-and-cut algorithm. All entries are
    rounded averages over three instances. CPU times are given in seconds.}\label{tab:random_prop}
  \scriptsize%
  \newcolumntype{P}{@{\extracolsep{\fill}}r@{\extracolsep{0ex}\;\;\;}r}%
  \newcolumntype{Q}{@{\extracolsep{\fill}}r@{\extracolsep{0ex}\;\;\;}r@{\;\;\;}r}%
  \begin{tabular*}{\textwidth}{@{}rrr@{\;\;}PPQ@{}}\toprule
             &        & & \multicolumn{2}{c}{basic} & \multicolumn{2}{c}{Iso Pruning} & \multicolumn{3}{c}{OF} \\
  $n$&   $m$  &   $q$ &         nsub &     cpu &          nsub &      cpu &         nsub &      cpu &     \#OF \\
\midrule
40 & 360 &  3   &           677  &           112  &            708  &          100  & \textbf{  516} & \textbf{  86} &              4   \\
40 & 360 &  6   &          1072  &            76  &            655  &           25  & \textbf{  157} & \textbf{  15} &             97   \\
40 & 360 &  9   & \textbf{    1} & \textbf{    0} & \textbf{     1} & \textbf{   0} & \textbf{    1} & \textbf{   0} &              0   \\
40 & 360 & 12   & \textbf{    1} & \textbf{    0} & \textbf{     1} & \textbf{   0} & \textbf{    1} & \textbf{   0} &              0   \\
\addlinespace
40 & 540 &  3   &           288  &           184  &            257  &          180  & \textbf{  219} & \textbf{ 146} &              4   \\
40 & 540 &  6   &         57606  &         13915  &          48786  &         7024  & \textbf{32347} & \textbf{5548} &           5750   \\
40 & 540 &  9   &         62053  &         14400  &         162871  &         7182  & \textbf{43434} & \textbf{4768} &          25709   \\
40 & 540 & 12   &         40598  &          6018  &           2187  &           69  & \textbf{  174} & \textbf{  31} &            166   \\
\addlinespace
40 & 720 &  3   &           488  &          1399  &            393  &         1325  & \textbf{  366} & \textbf{1080} &              5   \\
40 & 720 &  6   &          6888  &         11139  &           3507  &         3957  & \textbf{ 2563} & \textbf{3263} &            756   \\
40 & 720 &  9   &         21746  &         14400  &          20743  &         8220  & \textbf{12753} & \textbf{6820} &          10656   \\
40 & 720 & 12   &         24739  &         14400  &          68920  &         9164  & \textbf{21067} & \textbf{6209} &          23532   \\
 \bottomrule
\end{tabular*}
\end{table}

\begin{table}
  \caption{Results of the branch-and-cut algorithm. CPU times are given in seconds.}
  \label{tab:random_prop2}
  \scriptsize%
  \newcolumntype{P}{@{\extracolsep{\fill}}r@{\extracolsep{0ex}\;\;\;}r}%
  \newcolumntype{Q}{@{\extracolsep{\fill}}r@{\extracolsep{0ex}\;\;\;}r@{\;\;\;}r}%
  \begin{tabular*}{\textwidth}{@{}rrrr@{\;\;}PPQ@{}}
  \toprule
               &        &        & & \multicolumn{ 2 }{c}{basic} & \multicolumn{ 2 }{c}{Iso Pruning} & \multicolumn{ 3 }{c}{OF} \\
              graph & $n$ & $m$ & $q$ &    cpu &   nsub &    cpu &  nsub &   cpu & nsub & \#OF \\
  \midrule
 data\_2g\_10\_1001 & 100 & 200 & 3  & \textbf{   9} & \textbf{     7} &           11  &            17  & \textbf{  9} &           15  &    1 \\
 data\_2g\_10\_1001 & 100 & 200 & 5  &           84  &            674  & \textbf{  31} & \textbf{  148} &          48  &          346  &   98 \\
 data\_2g\_10\_1001 & 100 & 200 & 7  &           71  &             89  & \textbf{  30} &           233  &          52  & \textbf{  35} &   42 \\
    data\_2g\_6\_66 &  36 &  72 & 3  & \textbf{   1} & \textbf{    11} & \textbf{   1} &            27  & \textbf{  1} & \textbf{  11} &    3 \\
    data\_2g\_6\_66 &  36 &  72 & 5  &            5  &            483  &            5  &           667  & \textbf{  3} & \textbf{ 175} &   67 \\
    data\_2g\_6\_66 &  36 &  72 & 7  &            8  &            623  &            9  &          1259  & \textbf{  5} & \textbf{ 226} &  232 \\
  data\_2g\_7\_1034 &  49 &  98 & 3  & \textbf{   1} &             23  &            4  &            28  & \textbf{  1} & \textbf{  13} &    1 \\
  data\_2g\_7\_1034 &  49 &  98 & 5  &           10  &            613  & \textbf{   8} &           534  &          11  & \textbf{ 491} &   51 \\
  data\_2g\_7\_1034 &  49 &  98 & 7  &           50  &           3541  &           19  &          1229  & \textbf{ 11} & \textbf{ 359} &  110 \\
   data\_2g\_8\_648 &  64 & 128 & 3  & \textbf{   4} &             21  &            5  &            23  &           5  & \textbf{  15} &    2 \\
   data\_2g\_8\_648 &  64 & 128 & 5  &           15  &            383  & \textbf{   7} & \textbf{   10} &          12  &          159  &   18 \\
   data\_2g\_8\_648 &  64 & 128 & 7  & \textbf{   9} &             95  &           13  & \textbf{   12} &          72  &         1879  &  530 \\
  data\_2g\_9\_9211 &  81 & 162 & 3  & \textbf{   5} & \textbf{    48} &           22  &            77  &          11  &           61  &    1 \\
  data\_2g\_9\_9211 &  81 & 162 & 5  &           13  & \textbf{     1} & \textbf{  10} & \textbf{    1} &          15  & \textbf{   1} &    0 \\
  data\_2g\_9\_9211 &  81 & 162 & 7  & \textbf{  14} & \textbf{     1} &           17  & \textbf{    1} &          15  & \textbf{   1} &    0 \\
 \addlinespace
 data\_3g\_234\_234 &  24 &  60 & 3  & \textbf{   0} & \textbf{     5} & \textbf{   0} &            16  & \textbf{  0} & \textbf{   5} &    1 \\
 data\_3g\_234\_234 &  24 &  60 & 5  & \textbf{   0} & \textbf{     1} & \textbf{   0} & \textbf{    1} & \textbf{  0} & \textbf{   1} &    0 \\
 data\_3g\_234\_234 &  24 &  60 & 7  & \textbf{   0} & \textbf{     1} & \textbf{   0} & \textbf{    1} & \textbf{  0} & \textbf{   1} &    0 \\
 data\_3g\_244\_244 &  32 &  80 & 3  & \textbf{   1} &             61  & \textbf{   1} &           100  & \textbf{  1} & \textbf{  57} &    2 \\
 data\_3g\_244\_244 &  32 &  80 & 5  & \textbf{   0} & \textbf{     1} & \textbf{   0} & \textbf{    1} & \textbf{  0} & \textbf{   1} &    0 \\
 data\_3g\_244\_244 &  32 &  80 & 7  & \textbf{   1} & \textbf{     1} & \textbf{   1} & \textbf{    1} & \textbf{  1} & \textbf{   1} &    0 \\
 data\_3g\_333\_333 &  27 &  81 & 3  & \textbf{   0} & \textbf{     5} &            1  &            11  &           1  & \textbf{   5} &    1 \\
 data\_3g\_333\_333 &  27 &  81 & 5  & \textbf{   1} &             71  &           46  &           129  & \textbf{  1} & \textbf{  37} &   19 \\
 data\_3g\_333\_333 &  27 &  81 & 7  &            5  &            493  &         1984  &           463  & \textbf{  2} & \textbf{ 131} &  136 \\
 data\_3g\_334\_334 &  36 & 108 & 3  & \textbf{   1} &             41  &            2  &            56  & \textbf{  1} & \textbf{  29} &    3 \\
 data\_3g\_334\_334 &  36 & 108 & 5  & \textbf{   3} &            139  &            9  &           352  & \textbf{  3} & \textbf{  51} &   10 \\
 data\_3g\_334\_334 &  36 & 108 & 7  &           36  &           3697  &          795  &          1245  & \textbf{ 14} & \textbf{ 527} &  566 \\
 data\_3g\_344\_344 &  48 & 144 & 3  & \textbf{   2} &             47  &            5  &            64  & \textbf{  2} & \textbf{  29} &    3 \\
 data\_3g\_344\_344 &  48 & 144 & 5  &           24  &            887  &           35  &           830  & \textbf{ 17} & \textbf{ 300} &   53 \\
 data\_3g\_344\_344 &  48 & 144 & 7  &          182  &           7747  &          494  &          2889  & \textbf{ 50} & \textbf{1032} &  556 \\
 data\_3g\_444\_444 &  64 & 192 & 3  & \textbf{   5} & \textbf{    37} &           14  &            73  &           7  &           39  &    3 \\
 data\_3g\_444\_444 &  64 & 192 & 5  &          274  &           5491  &          184  &          3695  & \textbf{117} & \textbf{1377} &  272 \\
 data\_3g\_444\_444 &  64 & 192 & 7  &         5503  &         114011  &          528  &         13938  & \textbf{412} & \textbf{4815} & 2500 \\
\addlinespace
 \multicolumn{4}{c}{Total}           &         6338  &         139349  &         4292  &         28131  &         898  &        12225  &      \\
 \bottomrule
 \end{tabular*}
\end{table}

\begin{table}
  \scriptsize%
  \caption{Comparison between the CPLEX base variant (CPLEX-basic), CPLEX symmetry breaking (CPLEX-sym5), and OF (CPLEX-OF), averaged over three instances}
  \label{tab:cplex}
  \newcolumntype{P}{@{\extracolsep{\fill}}r@{\extracolsep{0ex}\;\;\;}r}%
  \begin{tabular*}{\textwidth}{@{}rrr@{\;\;}PPP@{}}\toprule
              &        & & \multicolumn{2}{c}{CPLEX-basic} & \multicolumn{2}{c}{CPLEX-sym5} & \multicolumn{2}{c}{CPLEX-OF} \\
$n$& $m$ &$q$ &          cpu &             nsub &          cpu &            nsub &          cpu &            nsub   \\
\midrule
24 & 190 & 2  &            0  &             837  & \textbf{  0} & \textbf{   758} &           0  &            774   \\
24 & 190 & 3  &            5  &            7534  &           4  &           5824  &           4  & \textbf{  5745}  \\
24 & 190 & 4  &           13  &           23445  &           8  &          13274  & \textbf{  8} & \textbf{ 12054}  \\
24 & 190 & 5  &           12  &           20284  &           6  &           9211  & \textbf{  4} & \textbf{  5819}  \\
24 & 190 & 6  &            7  &            9749  &           5  &           6610  & \textbf{  2} & \textbf{  2336}  \\
\addlinespace
24 & 250 & 2  & \textbf{   2} & \textbf{   2768} &           2  &           2806  &           3  &           3056   \\
24 & 250 & 3  &           42  &           51085  & \textbf{ 37} & \textbf{ 40281} &          47  &          42385   \\
24 & 250 & 4  &          263  &          341511  & \textbf{134} & \textbf{146861} &         139  &         156559   \\
24 & 250 & 5  &          540  &          749898  &         397  &         472858  & \textbf{191} & \textbf{220726}  \\
24 & 250 & 6  &          927  &         1224157  &         226  &         265439  & \textbf{154} & \textbf{171314}  \\
\addlinespace
30 & 200 & 2  &            1  & \textbf{    898} & \textbf{  1} &            990  &           1  &            953   \\
30 & 200 & 3  &            5  &            7134  & \textbf{  4} & \textbf{  5668} &           5  &           5899   \\
30 & 200 & 4  &            4  &            6450  &           3  &           3697  &           3  & \textbf{  3100}  \\
30 & 200 & 5  &            2  &            1981  & \textbf{  1} &            750  &           1  & \textbf{   706}  \\
30 & 200 & 6  &            0  &             359  &           0  &            180  &           0  & \textbf{   151}  \\
\addlinespace
30 & 233 & 2  &            1  &            1732  & \textbf{  1} & \textbf{  1663} &           2  &           1823   \\
30 & 233 & 3  &           15  &           18346  & \textbf{ 11} & \textbf{ 11936} &          15  &          13293   \\
30 & 233 & 4  &           20  &           29774  &          14  &          17736  & \textbf{ 12} & \textbf{ 14387}  \\
30 & 233 & 5  &           26  &           34350  &          11  &          13466  & \textbf{  9} & \textbf{  9014}  \\
30 & 233 & 6  &            7  &            7295  &           8  &           7180  & \textbf{  2} & \textbf{  1957}  \\
\addlinespace
30 & 266 & 2  & \textbf{   3} &            3177  &           3  & \textbf{  3120} &           3  &           3317   \\
30 & 266 & 3  &           49  &           46779  & \textbf{ 39} & \textbf{ 36642} &          46  &          36840   \\
30 & 266 & 4  &          119  &          144671  &          90  &          94227  & \textbf{ 77} & \textbf{ 75152}  \\
30 & 266 & 5  &          181  &          224103  &         110  &         120573  & \textbf{ 72} & \textbf{ 71303}  \\
30 & 266 & 6  &          263  &          267457  &         101  &          97577  & \textbf{ 54} & \textbf{ 39315}  \\
\addlinespace
30 & 300 & 2  & \textbf{   7} & \textbf{   6176} &           7  &           6495  &           8  &           7026   \\
30 & 300 & 3  &          204  &          186858  & \textbf{175} & \textbf{145598} &         222  &         157670   \\
30 & 300 & 4  &          738  &          766173  &         427  &         396841  & \textbf{413} & \textbf{378413}  \\
30 & 300 & 5  &         1283  &         1357270  &         456  &         441720  & \textbf{408} & \textbf{364906}  \\
30 & 300 & 6  &         1002  &          913069  &         362  &         326654  & \textbf{145} & \textbf{113401}  \\
\addlinespace
\multicolumn{3}{c}{Total} & 5741  &     6455320  &        2643  &        2696635  &        2051  &        1919394   \\\bottomrule
  \end{tabular*}

\end{table}

\begin{figure}[b]
  \centering
  \newlength{\plotheight}
  \setlength{\plotheight}{2.4cm}
  \newcommand{\mystyle}[1]{\tiny {#1}}
  \psfrag{250}[Bc][Bl]{\mystyle{250\,s}\quad}
  \psfrag{500}[Bc][Bl]{\mystyle{500\,s}\quad}
  \psfrag{750}[Bc][Bl]{\mystyle{750\,s}\quad}
  \psfrag{1000}[Bc][Bl]{\mystyle{1000\,s}\quad}
  \psfrag{ 0}{}
  \psfrag{ 3}[tc][Bc]{\mystyle{3}}
  \psfrag{ 6}[tc][Bc]{\mystyle{6}}
  \psfrag{ 9}[tc][Bc]{\mystyle{9}}
  \psfrag{ 12}[tc][Bc]{\mystyle{12}}
  \psfrag{ 13}{}
  \psfrag{3600}[Bc][Bl]{}
  \psfrag{7200}[Bc][Bl]{\mystyle{2\,h}\;}
  \psfrag{10800}[Bc][Bl]{}
  \psfrag{14400}[Bc][Bl]{\mystyle{4\,h}\;}
  \psfrag{50}[Bc][Bl]{\mystyle{50\,\%}\quad\;}
  \psfrag{100}[Bc][Bl]{\mystyle{100\,\%}\quad\;}
  \psfrag{150}[Bc][Bl]{\mystyle{150\,\%}\quad\;}
  \psfrag{10}[Bc][Bl]{}
  \psfrag{20}[Bc][Bl]{}
  \psfrag{30}[Bc][Bl]{\mystyle{30\,\%}\quad\;}
  \includegraphics[height=\plotheight]{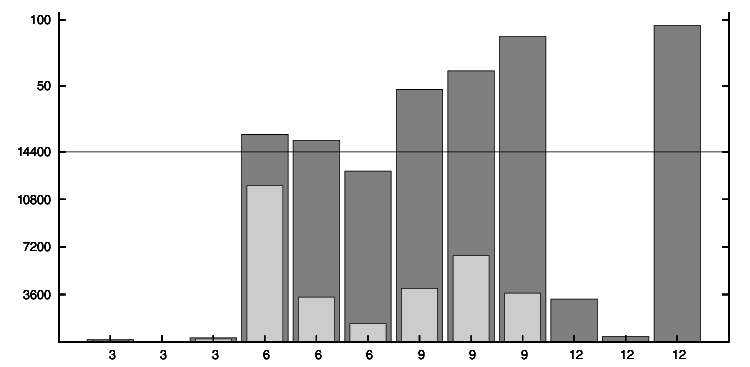}\hfill
  \psfrag{ 3}[tc][Bc]{\mystyle{3}}
  \psfrag{ 6}[tc][Bc]{\mystyle{6}}
  \psfrag{ 9}[tc][Bc]{\mystyle{9}}
  \psfrag{ 12}[tc][Bc]{\mystyle{12}}
  \includegraphics[height=\plotheight]{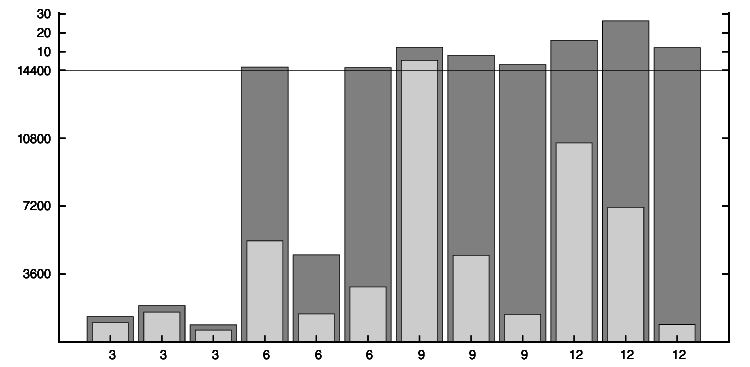}\hfill
  \psfrag{ 3}[tc][Bc]{\mystyle{3}}
  \psfrag{ 6}[tc][Bc]{\mystyle{6}}
  \psfrag{1800}[Bc][Bl]{\mystyle{1800\,s}\;}
  \psfrag{3600}[Bc][Bl]{}
  \psfrag{5400}[Bc][Bl]{\mystyle{5400\,s}\quad}
  \psfrag{7200}[Bc][Bl]{\mystyle{2\,h}}
  \psfrag{10800}[Bc][Bl]{}
  \psfrag{14400}[Bc][Bl]{\mystyle{4\,h}}
  \psfrag{50}[Bc][Bl]{\mystyle{50\,\%}\quad\;}
  \psfrag{100}[Bc][Bl]{\mystyle{100\,\%}\quad\;}
%
  \caption{Computation times/gaps for the basic version (dark gray) and the
    version with orbitopal fixing (light gray). Left: instances
    with $n=40$, $m = 540$. Right: instances for $n = 40$, $m = 720$.
    The number of partitions~$q$ is
    indicated on the $x$-axis. Values above 4 hours indicate the gap in
    percent.}
  \label{fig:basic_vs_OF}
\end{figure}

In a first experiment we tested the speedup that can be obtained by
performing orbitopal fixing. For this we compare the variant (\emph{basic})
without symmetry breaking (except for the zero-fixing of the upper right
$x$-variables) and the version in which we use orbitopal fixing
(\emph{OF}); see Table~\ref{tab:random_prop} for the results. Columns
``nsub'' give the number of nodes in the branch-and-bound tree and
``\#OF'' the number of fixings within OF. The results show that orbitopal
fixing is clearly superior (OF winners: 30, basic winners: 0), see also
Figure~\ref{fig:basic_vs_OF}.

Table~\ref{tab:random_prop} shows that the sparse instances are extremely
easy, the instances with $m = 540$ are quite easy, while the dense
instances are hard. A situation that often occurs for small~$m$ and
large~$q$ is that the optimal solution is 0, and hence no work has to be
done. For $m = 720$, the hardest instances arise when $q = 9$. It seems
that for~$q=3$ the small number of variables helps, while for $q = 12$ the
small objective function values help. Of course, symmetry breaking methods
become more important when~$q$ gets larger.

In addition, we report results for instances used in~\cite{GahAL09} in
Table~\ref{tab:random_prop2}. These instances arise from grid graphs. Due
to their sparsity these instances are already approachable for the basic
variant, but symmetry breaking yields huge performance gains here, too. As
sparsity is particularly exploited in our formulation, the running times
turn out to be much smaller than in the SDP based approach of
\cite{GahAL09}.

To compare orbitopal fixing to the isomorphism pruning approach of Margot,
we implemented the basic variant with fixed canonical variable order, the
\emph{ranked branching rule} (see \cite{Mar03}), as well as the variant
with free branching decisions, each adapted to the special symmetry we
exploit, which simplifies Margot's algorithm significantly. We decided to
use the canonical order variant, as it gave the best results. Other than
that, the same implementation and settings were used. It can be seen from
Table~\ref{tab:random_prop} (columns \emph{Iso Pruning}) that isomorphism
pruning is inferior to both orbitopal fixing (OF winners: 30, isomorphism
pruning winners: 0) and shifted column inequalities (30:0), but is still a
big improvement over the basic variant (28:2). Table~\ref{tab:random_prop2}
yields a similar conclusion: orbitopal fixing outperforms the basic
variant, as well as isomorphism pruning both in terms of cpu time 
(OF winners: 22, isomorphism pruning winners: 12, basic winners: 19),
and in terms of branch-and-bound nodes (OF : isomorphism pruning : basic = 27 : 9 : 12).
 We additionally report the total values over all instances in Table~\ref{tab:random_prop2}.
 Note   that there
are many easy instances in this instance set.
All in all, also in these experiments orbitopal fixing turns out to be superior and 
isomorphism
pruning still shows advantages over the basic variant.

We did not implement orbital branching, since this method when using global
symmetry only is very similar to isomorphism pruning in our context as
pointed out in Section~\ref{sec:Margot}. It should, however, be noted that,
in contrast to orbitopal fixing, both isomorphism pruning as well as
orbital branching could exploit symmetries of the instance graphs, too.
However, no nontrivial graph automorphisms in our test instances were found
by {\ttfamily nauty}, see~\cite{McK81}.

In a second experiment, we investigated the symmetry breaking capabilities
built into CPLEX. We suspect that it breaks symmetry within the tree,
but no detailed information was available. We first ran CPLEX
12.1 on the
IP formulation stated in Sect.~\ref{sec:introduction}. In one variant, we
fixed variables~$x_{ij}$ with~$j > i$ to zero, but turned symmetry breaking
off.  In a second variant, we turned symmetry breaking on and did not fix
variables to zero (otherwise CPLEX seems not to recognize the symmetry).
The level was set to most aggressive (5), although the default setting
yields the same results. The symmetry breaking variant turned out to be
effective: it was always faster than the basic version without symmetry
breaking. The black box use of
CPLEX was always inferior to our code. However, of course this
is partially due to our use of specialized cutting planes.

We then compared the built-in symmetry treatment of CPLEX to orbitopal
fixing. We implemented orbitopal fixing as a branching callback in CPLEX
that passed all fixings with the branching decision. In order to obtain a
fair comparison, all advanced features like preprocessing, primal
heuristics, and cuts were turned off.  Moreover, we did not provide an
optimal solution.  The branching was performed as a first--index branching
along the $x$--variables. The results are shown in Table~\ref{tab:cplex}
for several different graph densities. The comparison is performed between CPLEX without symmetry
handling (\emph{CPLEX-basic}), with aggressive symmetry handling
(\emph{CPLEX-sym5}), and with orbitopal fixing (\emph{CPLEX-OF}).
For the reasons  discussed above, we only fixed the upper triangle to zero in variants \emph{CPLEX-basic} and \emph{CPLEX-OF} .
It can be seen that in these experiments
orbitopal fixing performs slightly better than symmetry breaking in CPLEX.
Note that the comparison is skewed in favor of CPLEX's own symmetry
handling: orbitopal fixing relies on the ability to identify the fixed
variables in every node of the branch-and-bound tree.  However, apart from
the branching decisions no fixings are reported by the CPLEX-API. In contrast,
SCIP is much better suited for orbitopal fixing, since
in
SCIP the strengthening of variable bounds, called
propagation, is an essential concept. Besides this, SCIP can also perform conflict
analysis, which makes use of the information collected via propagation,
see~\cite{Ach07b}.

In another experiment, we turned off orbitopal fixing and separated shifted
column inequalities in every node of the tree. The results on the original
testset of random instances are that the OF-version is slightly better than
the SCI variant (OF winners: 19, SCI winners: 11), but the results are
quite close (OF average time: 2330 seconds, SCI average time: 2288
seconds). Although by Part~2 of Theorem~\ref{thm:seqfixOrbi}, orbitopal
fixing is as strong as fixing with SCIs (with the same branching
decisions), the LPs get harder and the process slows down a bit. On the
other hand, the SCIs are already active in the root node, which in general yields a
better root bound. This may result in fewer branch-and-bound nodes due to
potentially more fixings in the root node.

\section{Concluding Remarks}

The main contribution of this paper is the development of an algorithm that
handles orbitopal symmetry for binary programs with assignment structure by
fixing values in partial solutions to exclude symmetric branches from
exploration.  The algorithm is proven to be optimal in the sense that as
many such fixings, as early as possible are made. Moreover, it is shown
that the algorithm can be implemented to run in linear time. The considered
assignment structure occurs frequently in standard IP formulations.

The effectiveness of our approach is demonstrated by applying our
algorithm to the Graph Partitioning problem, where it is compared
with other known methods to handle symmetry.

In the future, extensions in several directions might be possible:
other group actions, different restrictions on the number of $1$'s in
each row, symmetries acting on both rows and columns.
Moreover, more extensive
studies for other applications are desirable to further explore the
practical implications of our symmetry handling approach.

\section*{Acknowledgments}

We are grateful to Andreas Loos for discussions on the complexity of
optimizing linear functions over covering orbitopes that lead to the basic
idea for the proof of Theorem~5. We also thank Bissan Ghaddar for providing
us with the instances reported on in Table~\ref{tab:random_prop2} and an
anonymous referee for comments that helped to improve the computational
results of this paper.

\bibliographystyle{siam}

\end{document}